\documentclass[11pt,leqno]{article}

\usepackage[all]{xy}				
\usepackage{amsthm}				
\usepackage{amsfonts}				
\usepackage{bussproofs}		     		
\usepackage{cancel}				
\usepackage{fancyhdr}				
\usepackage{geometry}				
\usepackage{graphicx}				
\usepackage{hyperref}				
\usepackage[utf8]{inputenc}			
\usepackage{lineno}				
\usepackage{mathrsfs}				
\usepackage{MnSymbol}			    	
\usepackage{tikz}				
\usetikzlibrary{tikzmark}
\usepackage{url}				
\usepackage{xcolor}				

\usepackage[numbers]{natbib}

\theoremstyle{definition}      		
\newtheorem{df}{Definition}    		
\newtheorem{cor}[df]{Corollary}		
\newtheorem{lm}[df]{Lemma}	  	
\newtheorem{prop}[df]{Proposition}	
\newtheorem{rmrk}[df]{Remark}		
\newtheorem{tr}[df]{Theorem}	 	

	\newif\ifSuppressEndOfProva\SuppressEndOfProvafalse

\usepackage[shortlabels]{enumitem}

\newcommand{\noi}{\noindent}
\newcommand{\seta}{\renewcommand{\arraystretch}{0.95}}
\newcommand{\oid}{\ensuremath{\Leftrightline}}
\newcommand{\noid}{\ensuremath{\nLeftrightline}}
\newcommand{\D}{\ensuremath{\mathcal{D}}}
\newcommand{\I}{\ensuremath{\mathcal{I}}}
\newcommand{\sig}{\ensuremath{\mathcal{S}}}
\newcommand{\ov}{\ensuremath{\overline}}
\newcommand{\wh}{\ensuremath{\widehat}}

\label{xLOGICS}{}

\newcommand{\letf}{$LET_{F}$}
\newcommand{\qletf}{$QLET_{F}$}

\newcommand{\fde}{$FDE$}
\newcommand{\qfde}{$QFDE$}

\newcommand{\qmbc}{$QmbC$}

\newcommand{\nel}{\textit{N4}}
\newcommand{\lets}{\textit{LET}s}
\newcommand{\qlp}{\ensuremath{QLP}}
\newcommand{\qk}{\ensuremath{\textit{QK3}}}
\newcommand{\qcl}{\ensuremath{QCL}}

\label{xCOLORS} %

\newcommand{\black}{\color{black}}

\label{xSYMBOLS}

\newcommand{\cons}{\ensuremath{{\circ}}}
\newcommand{\con}{\ensuremath{{\circ}}}
\newcommand{\incon}{\ensuremath{{\bullet}}}

 \newcommand{\A}{\ensuremath{\mathfrak{A}}}

\label{xGENERAL COMMANDS}%

\newcommand{\mh}{\noindent}

\newcommand{\setl}{\setlength\itemsep{-0.2em}}

\newcommand{\bqu}{\begin{quote}} 
\newcommand{\equ}{\end{quote}}
\newcommand{\enr}{\begin{enumerate}[label={(\arabic*)}, resume]}
\newcommand{\eenr}{\end{enumerate}}

\label{xPORTUGUES}

\begin{document}


\title{\textbf{
Valuation semantics for first-order  logics of evidence and truth (and some related logics)
}
\thanks{
We would like to thank Mart\'in Figallo and Andrea Loparic for   valuable discussions that helped to shape  some 
of the ideas presented in  this text. }
}

\author{H. Antunes, A. Rodrigues, W. Carnielli, M. E. Coniglio}

\maketitle

\begin{abstract} \noi This paper  introduces  the logic \qletf, a quantified extension of the logic of
evidence and truth \letf,
together  with a corresponding sound and complete first-order non-deterministic valuation semantics.
\letf\ is a paraconsistent and paracomplete sentential logic that extends the logic of first-degree
entailment (\fde) with a classicality operator \cons\ and  a non-classicality operator \incon,
dual to each other: while $\con A$ entails that $A$ behaves classically, $\incon A$ follows from $A$'s
violating some classically valid inferences.  The semantics of \qletf\ combines structures that interpret
negated predicates in terms of  anti-extensions with   first-order non-deterministic valuations, and
completeness is  obtained through a generalization of Henkin's method.  By providing sound and complete
semantics for first-order extensions of  \fde, \textit{K3}, and \textit{LP}, we show how these tools,
which we call here the method of \textit{{anti-extensions + valuations}}, can be naturally applied  to a
number of non-classical logics.
\end{abstract}

\section*{Introduction}\label{sec:introduction}

The main aim of this paper is to introduce the logic \qletf, a quantified extension of the logic of
evidence and truth \letf, introduced in Rodrigues,  Bueno-Soler, \& Carnielli \cite{letf}, together  with
a corresponding sound and complete first-order valuation semantics.  The latter  are a development of the
non-deterministic  semantics for sentential logics investigated by Loparic et al.
\cite{costa.alves,loparic1986,loparic.alves,costa_loparic} from the 1970s onward in order to provide
adequate semantics for some non-classical logics. \letf\ is a paraconsistent and paracomplete sentential
logic that extends the logic of first-degree entailment (\fde), also known as Belnap-Dunn 4-valued logic,
with a classicality operator \cons\ and a non-classicality operator \incon, dual to \con\ in the sense
that  $\con A$ entails that $A$ behaves classically, and $\incon A$ follows from $A$'s violating some
classically valid inferences.  A sound and complete valuation semantics for \letf\ was presented in
\cite{letf}, and a Kripke-style semantics in \cite{axioms}.

Logics of evidence and truth (\lets) have been conceived to formalize the deductive behavior of positive
and negative evidence, which can be either conclusive or non-conclusive\footnote{For a more detailed
discussion of the notion of evidence, see \cite[Sect.~2]{barrio.re}.}.  \lets\ consider that conclusive
evidence behaves classically, and so is subjected to classical logic. Non-conclusive evidence, on the
other hand, may be incomplete or contradictory and, in the case of \letf, is subjected to \fde.
According to the intended interpretation in terms of evidence and truth, a pair of contradictory
sentences $A$ and $\neg A$ expresses conflicting non-conclusive evidence for $A$, and $\con A$ expresses
that there is conclusive evidence for either the truth or the falsity of $A$.

The semantics of \qletf\ combines structures that interpret  negated predicates in terms of
anti-extensions \cite[e.g.][]{henr.tese, priest2002} with first-order non-deterministic valuations
\cite[e.g.][]{carco.book,qmbc},  and completeness is obtained by a Henkin-style proof.  These tools,
which for convenience we call here the method of \textit{{anti-extensions + valuations}}, are required
for handling the non-deterministic character of \qletf, and can be naturally applied  to a number of
non-classical logics.

It is well known that  \fde\ can be interpreted as an information-based logic
\cite[e.g.][]{belnap1977.how,dunn76,dunn2019}. \label{page.info} It has been argued in \cite{axioms,
letf} that \letf\ can be interpreted along the same lines, and this interpretation can be naturally
extended to \qletf. The latter may be seen as representing databases that contain not only contradictory
and incomplete information but also reliable (true) information.\footnote{The notion of evidence can be
explained based on the notion of information as \textit{meaningful data} -- see~e.g.~Fetzer
\cite[][]{fetzer2004a}.  In line with Fetzer,  Dunn in \cite[p.~589]{dunn2008} thinks of information as a
pure semantic content that may be false and does not depend on the belief of an agent. Evidence is thus
just information that comes together with a justification that may fail to justify that information  it
is intended to justify \cite[cf.][Sect.~2.3]{barrio.re}.}
In this case, $\cons A$ means that the information about $A$ is reliable, while $\incon A$ means that
there is no reliable information about $A$.

A quantified version of \fde\ is obtained as a fragment of the logic \qletf. The former can, in turn, be
extended to obtain quantified versions of the well-known Kleene's  \textit{K3} \cite{kleene} and the
logic of paradox \textit{LP} \cite{priest.lp}
by adding rules corresponding respectively to the principles of explosion and excluded middle. We will
show how the method of {anti-extensions + valuations }  can be adapted to provide sound and complete
semantics for these logics, called here, respectively, \qfde, \qk\ and \qlp,   and also for classical
logic, that is a sort of limiting case where anti-extensions and non-deterministic valuations are not
necessary.

The remainder of this paper is structured as follows. In Section \ref{sec:qletf} we present both a
natural deduction system and a corresponding valuation semantics for \qletf. Section
\ref{sec:soundnessandcompleteness} contains detailed proofs of soundness and completeness results for
\qletf, along with other metatheoretical results, such as compactness and a few versions of the
L\"owenheim-Skolem theorem. In Section 3 we provide sound and complete semantics for the first-order
versions of \fde, \textit{K3}, and \textit{LP}, which are particular cases of the  method of
{anti-extensions + valuations}.  Section 4 wraps up the text with some historical remarks about   
valuation semantics and their generality.

\section{The logic \qletf}\label{sec:qletf}

The logical vocabulary of \qletf\ is composed by the unary connectives $\neg$, $\con$, $\incon$, the
binary connectives $\land$ and $\lor$, the quantifiers $\forall$ and $\exists$, the identity symbol
$\oid$, the individual variables from $\mathcal{V} = \{v_{i}: i \in \mathbb{N}\}$, and parentheses. From
now on we shall specify the non-logical vocabulary of a first-order language by means of its
\textit{signature}, which is a pair $\sig = \langle \mathcal{C}, \mathcal{P} \rangle$ such that
$\mathcal{C}$ is an infinite set of individual constants and $\mathcal{P}$ is a set of predicate letters.
Each element $P$ of $\mathcal{P}$ is assumed to have a corresponding finite arity and $\oid \ \in
\mathcal{P}$ is a binary predicate.  Given a signature \sig, its cardinality is the cardinality of the
set $\mathcal{C} \cup \mathcal{P}$.

Henceforth, we implicitly assume the usual definitions of such syntactic notions as \emph{term},
\emph{formula}, \emph{bound/free occurrence of a variable}, \emph{sentence} etc.  -- but with the
proviso that formulas with void quantifiers are not allowed.   Given a signature \sig, we shall
denote the set of terms generated by \sig\ by $Term(\sig)$. Likewise, the set of formulas and the set of
sentences generated by \sig\ will be denoted by $Form(\sig)$ and $Sent(\sig)$\footnote{Hereafter $x$,
$x_{1}$, $x_{2}$ will be used as metavariables ranging over $\mathcal{V}$, $c$, $c_{1}$, $c_{2}$,$\dots$
as metavariables ranging over $\mathcal{C}$, $t$, $t_{1}$, $t_{2}$, $\dots$ as metavariables ranging over
$Term(\sig)$, and $A$, $B$, $C$, $\dots$ as metavariables ranging over $Form(\sig)$. Given $t, t_{1},
t_{2} \in Term(\sig)$, we will use the notation $t(t_{2}/t_{1})$ to denote the result of replacing every
occurrence of $t_{1}$ in $t$ (if any) by $t_{2}$. Similarly, $A(t/x)$ will denote the formula that
results by replacing every free occurrence of $x$ in $A$ by $t$.}.

For the sake of simplicity, the deductive systems and the formal semantics of the logics discussed below
will be formulated exclusively in terms of sentences. This is the reason why we've assumed right from the
outset that languages must always have an infinite stock of individual constants -- for otherwise we
could be prevented from applying some of the quantifier rules   due to the lack of enough constants.
However, none of the following definitions and results depend essentially on this decision (see Remark
\ref{rmrk:openformulas}).

	\begin{df}\label{df:nds} Let \sig\ be a signature, $c \in \mathcal{C}$, and $A, B, C \in
	Sent(\sig)$. The logic \qletf\ is defined over \sig\ by the following natural deduction rules:

		{\small

			\begin{center}
				\bottomAlignProof
					\AxiomC{$A$}
						\AxiomC{$B$}
					\RightLabel{$\land I$}
					\BinaryInfC{$A \land B$}
				\DisplayProof
			\qquad
				\bottomAlignProof
						\AxiomC{$A \land B$}
					\RightLabel{$\land E$}
					\UnaryInfC{$A$}
				\DisplayProof
			$\!\!$
				\bottomAlignProof
						\AxiomC{$A \land B$}
					\UnaryInfC{$B$}
				\DisplayProof

			\end{center}

			\vspace{0.25cm}

			\begin{center}

				\bottomAlignProof
						\AxiomC{$A$}
					\RightLabel{$\lor I$}
					\UnaryInfC{$A \lor B$}
				\DisplayProof
			$\!\!$
				\bottomAlignProof
						\AxiomC{$B$}
					\UnaryInfC{$A \lor B$}
				\DisplayProof
			\qquad
				\bottomAlignProof
						\AxiomC{$A \lor B$}
							\AxiomC{$[A]$} \noLine
						\UnaryInfC{$\vdots$} \noLine
						\UnaryInfC{$C$}
							\AxiomC{$[B]$} \noLine
						\UnaryInfC{$\vdots$} \noLine
						\UnaryInfC{$C$}
					\RightLabel{$\lor E$}
					\TrinaryInfC{$C$}
				\DisplayProof

			\end{center}

			\vspace{0.25cm}

			\begin{center}

				\bottomAlignProof
						\AxiomC{$\neg A$}
					\RightLabel{$\neg \land I$}
					\UnaryInfC{$\neg(A \land B)$}
				\DisplayProof
			$\!\!$
				\bottomAlignProof
						\AxiomC{$\neg B$}
					\UnaryInfC{$\neg(A \land B)$}
				\DisplayProof
			\qquad
				\bottomAlignProof
						\AxiomC{$\neg(A \land B)$}
							\AxiomC{$[\neg A]$} \noLine
						\UnaryInfC{$\vdots$} \noLine
						\UnaryInfC{$C$}
							\AxiomC{$[\neg B]$} \noLine
						\UnaryInfC{$\vdots$} \noLine
						\UnaryInfC{$C$}
					\RightLabel{$\neg \land E$}
					\TrinaryInfC{$C$}
				\DisplayProof

			\end{center}

			\vspace{0.25cm}

			\begin{center}

				\bottomAlignProof
						\AxiomC{$\neg A$}
						\AxiomC{$\neg B$}
					\RightLabel{$\neg \lor I$}
					\BinaryInfC{$\neg(A \lor B)$}
				\DisplayProof
			\qquad
				\bottomAlignProof
						\AxiomC{$\neg(A \lor B)$}
					\RightLabel{$\neg \lor E$}
					\UnaryInfC{$\neg A$}
				\DisplayProof
			$\!\!$
				\bottomAlignProof
						\AxiomC{$\neg(A \lor B)$}
					\UnaryInfC{$\neg B$}
				\DisplayProof

			\end{center}

			\vspace{0.25cm}

			\begin{center}

				\bottomAlignProof
						\AxiomC{$A$}
					\RightLabel{$DN$}
					\UnaryInfC{$\neg \neg A$}
				\DisplayProof
			$\!\!$
				\bottomAlignProof
						\AxiomC{$\neg \neg A$}
					\UnaryInfC{$A$}
				\DisplayProof

			\end{center}

			\vspace{0.25cm}

			\begin{center}

				\bottomAlignProof
						\AxiomC{$\con A$}
						\AxiomC{$A$}
						\AxiomC{$\neg A$}
					\RightLabel{$EXP^{\con}$}
					\TrinaryInfC{$B$}
				\DisplayProof
			\qquad
				\bottomAlignProof
						\AxiomC{$\con A$}
					\RightLabel{$PEM^{\con}$}
					\UnaryInfC{$A \lor \neg A$}
				\DisplayProof

			\end{center}

			\vspace{0.25cm}

			\begin{center}

				\bottomAlignProof
						\AxiomC{$\con A$}
						\AxiomC{$\incon A$}
					\RightLabel{$Cons$}
					\BinaryInfC{$B$}
				\DisplayProof
			\qquad
				\bottomAlignProof
						\AxiomC{}
					\RightLabel{$Comp$}
					\UnaryInfC{$\con A \lor \incon A$}
				\DisplayProof

			\end{center}

			\begin{center}

				\bottomAlignProof
						\AxiomC{$B \lor A(c/x)$}
					\RightLabel{$\forall I$}
					\UnaryInfC{$B \lor \forall xA$}
				\DisplayProof
			\qquad
				\bottomAlignProof
						\AxiomC{$\forall xA$}
					\RightLabel{$\forall E$}
					\UnaryInfC{$A(c/x)$}
				\DisplayProof
			\qquad
				\bottomAlignProof
						\AxiomC{$A(c/x)$}
					\RightLabel{$\exists I$}
					\UnaryInfC{$\exists xA$}
				\DisplayProof
			\qquad
				\bottomAlignProof
						\AxiomC{$\exists xA$}
							\AxiomC{$[A(c/x)]$} \noLine
						\UnaryInfC{$\vdots$} \noLine
						\UnaryInfC{$C$}
					\RightLabel{$\exists E$}
					\BinaryInfC{$C$}
				\DisplayProof

			\end{center}

			\vspace{0.25cm}

			\begin{center}

				\bottomAlignProof
						\AxiomC{$\neg A(c/x)$}
					\RightLabel{$\neg \forall I$}
					\UnaryInfC{$\neg \forall xA$}
				\DisplayProof
				\qquad
				\bottomAlignProof
						\AxiomC{$\neg \forall xA$}
							\AxiomC{$[\neg A(c/x)]$} \noLine
						\UnaryInfC{$\vdots$} \noLine
						\UnaryInfC{$C$}
					\RightLabel{$\neg \forall E$}
					\BinaryInfC{$C$}
				\DisplayProof
			\qquad
				\bottomAlignProof
						\AxiomC{$\neg A(c/x)$}
					\RightLabel{$\neg \exists I$}
					\UnaryInfC{$\neg \exists xA$}
				\DisplayProof
				\qquad
				\bottomAlignProof
						\AxiomC{$\neg \exists xA$}
					\RightLabel{$\neg \exists E$}
					\UnaryInfC{$\neg A(c/x)$}
				\DisplayProof

			\end{center}

			\vspace{0.25cm}

			\begin{center}

				\bottomAlignProof
						\AxiomC{}
					\RightLabel{$\oid I$}
					\UnaryInfC{$c \oid c$}
				\DisplayProof
				\qquad
				\bottomAlignProof
						\AxiomC{$c_1 \oid c_2$}
						\AxiomC{$A(c_1/x)$}
					\RightLabel{$\oid E$}
					\BinaryInfC{$A(c_2/x)$}
				\DisplayProof
			\qquad
				\bottomAlignProof
						\AxiomC{$A$}
					\RightLabel{$AV$}
					\UnaryInfC{$A'$}
				\DisplayProof

			\end{center}
		}

	\noi In $\forall I$, $c$ must not occur in $A$ or $B$, nor in any hypothesis on which $B \lor
	A(c/x)$ depends; and in $\neg \exists I$, $c$ must not occur in $A$ nor in any hypothesis on
	which $\neg A(c/x)$ depends. In $\exists E$ and $\neg \forall E$, $c$ must occur neither in $A$
	or $C$, nor in any hypothesis on which $C$ depends, except $A(c/x)$ ($\neg A(c/x)$). Finally, in
	$AV$, $A'$ denotes any \emph{alphabetic variant of} $A$\footnote{A formula is an \emph{alphabetic
	variant} of another if they only differ in (some of) their bound variables. See \cite[pp.
	126-7]{Enderton2001}  for details.}.

	\end{df}

	\begin{prop}\label{prop:generalization} The usual universal generalization rule:

		\begin{center}

			\bottomAlignProof
					\AxiomC{$A(c/x)$}
				\UnaryInfC{$\forall xA$}
				\RightLabel{$\forall I'$}
			\DisplayProof

		\end{center}

	\noi (where $c$ occurs neither in $A$ nor in any hypothesis on which $A(c/x)$ depends) can be
	derived in \qletf.
	\end{prop}

	\begin{pr} It suffices to consider the following derivation:

		{\small \begin{center}

			\bottomAlignProof
							\AxiomC{$A(c/x)$}
						\RightLabel{$\lor I$}
						\UnaryInfC{$\forall xA \lor A(c/x)$}
					\RightLabel{$\forall I$}
					\UnaryInfC{$\forall xA \lor \forall xA$}
					\AxiomC{$[\forall xA]_{1}$}
					\AxiomC{$[\forall xA]_{1}$}
				\RightLabel{$\lor E_{1}$}
				\TrinaryInfC{$\forall xA$}
			\DisplayProof

		\end{center}
		}
	\end{pr}

Given a signature \sig\ and $\Gamma \cup \{A\} \subseteq Sent(\sig)$, the definition of a deduction of
$A$ from $\Gamma$ in \qletf\ is the usual one \cite[see e.g.][Ch. 2]{Troelstra&vanDalen1988}.  It
suffices to say here that a derivation $\Theta$ is a tree of labeled sentences in which each node either
is an element of the set of premises $\Gamma$ or results from preceeding nodes by the application of one
of the rules above, and whose bottomost sentence is the \textit{conclusion} of $\Theta$. We shall use the
notation $\Gamma \vdash_{\sig} A$ to express that there exists a derivation in \qletf\ from the premises
in $\Gamma$ and whose conclusion is $A$, omitting the subscript if there is no risk of confusion.

	\begin{prop}\label{prop:alternativerules} Consider the following eight rules:

		{\small \begin{center}
			\bottomAlignProof
					\AxiomC{$\forall x \neg A$}
				\UnaryInfC{$\neg \exists x A$}
			\DisplayProof
		\qquad
			\bottomAlignProof
					\AxiomC{$\neg \exists x A$}
				\UnaryInfC{$\forall x \neg A$}
			\DisplayProof
		\qquad
			\bottomAlignProof
					\AxiomC{$\exists x \neg A$}
				\UnaryInfC{$\neg \forall xA$}
			\DisplayProof
		\qquad
			\bottomAlignProof
					\AxiomC{$\neg \forall xA$}
				\UnaryInfC{$\exists x \neg A$}
			\DisplayProof
		\end{center}

		\vspace{0.25cm}

		\begin{center}
			\bottomAlignProof
					\AxiomC{$\forall xA$}
				\UnaryInfC{$\neg \exists x \neg A$}
			\DisplayProof
		\qquad
			\bottomAlignProof
					\AxiomC{$\neg \exists x \neg A$}
				\UnaryInfC{$\forall xA$}
			\DisplayProof
			\qquad
			\bottomAlignProof
					\AxiomC{$\exists xA$}
				\UnaryInfC{$\neg \forall x \neg A$}
			\DisplayProof
		\qquad
			\bottomAlignProof
					\AxiomC{$\neg \forall x \neg A$}
				\UnaryInfC{$\exists xA$}
			\DisplayProof
		\end{center}
		}

		\begin{itemize}\setl

		\item[1.] Each one of these rules can be derived in \qletf;

		\item[2.] The first four rules, together with $\forall I$, $\forall E$, $\exists I$, and
		$\exists E$, are sufficient for deriving $\neg \forall I$, $\neg \forall E$, $\neg
		\exists I$, and $\neg \exists E$.

		\end{itemize}
	\end{prop}

	\begin{pr} Left to the reader.
	\end{pr}

\subsection{First-order valuation semantics for \qletf}

	\begin{df}\label{df:structure} Let \sig\ be a signature. An \sig-\emph{structure} \A\ is a pair
	$\langle \D, \I \rangle$ such that \D\ is a non-empty set (the \textit{domain} of \A) and \I\ is
	an \textit{interpretation function} such that:

		\begin{itemize}\setl

		\item[1.] For every constant $c \in \mathcal{C}$, $\I(c) \in \D$;

		\item[2.] For every $n$-ary predicate $P \in \mathcal{P}$, $\I(P)$ is a pair $\langle
		P_{+}^{\A}, P_{-}^{\A} \rangle$ such that $P_{+}^{\A} \cup P_{-}^{\A} \subseteq \D^{n}$.

		\item[3.] $\oid_{+}^{\A} \ = \{\langle a, a \rangle: a \in \D\}$.

		\end{itemize}

	\noi Given an $\sig$-structure $\A = \langle \D, \I \rangle$, we shall write $c^{\A}$ and
	$P^{\A}$ instead of respectively $\I(c)$ and $\I(P)$.
	\end{df}

According to the definition above, individual constants are interpreted as elements of the domain \D\ of
\A, while predicate letters are interpreted as pairs of relations over \D: each predicate letter $P$,
including $\oid$, is assigned both an \emph{extension}, $P_{+}^{\A}$, and an \emph{anti-extension},
$P_{-}^{\A}$, where $P_{+}^{\A}$ and $P_{-}^{\A}$   are intended to express, respectively, the presence
of positive and negative evidence (or information) for the atomic sentences of the relevant language.


Notice that, given an $n$-ary predicate letter $P$, although $P_{+}^{\A} \cup P_{-}^{\A}$ must be a
subset of $\D^{n}$, there are no constraints to the effect that $P_{+}^{\A} \cup P_{-}^{\A} = \D^{n}$,
nor to the effect that $P_{+}^{\A} \cap P_{-}^{\A} = \emptyset$.  As it will become clear below, this
means that it is not required that exactly one of $P(c_{1},\dots,c_{n})$ and $\neg P(c_{1},\dots,c_{n})$
receive a designated value, for every constants $c_{1},\dots,c_{n}$: it may be that neither
$P(c_{1},\dots,c_{n})$ nor $\neg P(c_{1},\dots,c_{n})$ holds in \A, or that both $P(c_{1},\dots,c_{n})$
and $\neg P(c_{1},\dots,c_{n})$ do.  Notice further that albeit $\oid$ is also interpreted as a pair of
relations, as any other predicate letter, its extension, $\oid_{+}^{\A}$, must be the identity relation
on \D, which is meant to ensure that $\oid$ satisfy the most basic properties of identity.  Nonetheless,
nothing prevents there existing some $a \in \D$ such that $\langle a, a \rangle \in \ \oid_{-}^{\A}$, in
which case both $c \oid c$ and $\neg(c \oid c)$ will receive a designated value, for some $c$ (we return
to the interpretation of identity   in Remark \ref{rmk:identity}).

	\begin{df}\label{df:diagram} Let $\sig = \langle \mathcal{C}, \mathcal{P} \rangle$ be a signature
	and let $\A$ be an \sig-structure.  The \emph{diagram signature} $\sig_{\A}$ of $\A$ is the pair
	$\langle \mathcal{C}_{\A}, \mathcal{P} \rangle$ such that $\mathcal{C}_{\A} = \mathcal{C} \cup
	\{\ov{a}: a \in \D\}$; that is, $\sig_{\A}$ is the signature that results from $\sig$ by
	introducing a new individual constant $\ov{a}$ for each element $a$ of the domain.  The language
	generated by $\sig_{\A}$ will be called \emph{the diagram language of} $\A$, and we shall use the
	notation $\wh{\A}$ to denote the $\sig_{\A}$-structure that is just like \A\ except that
	$\ov{a}^{\wh{\A}} = a$, for every $a \in \D$.  \end{df}

	\begin{df}\label{df:valuation} Let $\sig = \langle \mathcal{C}, \mathcal{P}\rangle$ be a signature and
	$\A$ be an \sig-structure.  A mapping $v: Sent(\sig_{\A}) \longrightarrow \{0,1\}$ is an
	$\A$-\emph{valuation} if it satisfies the following conditions:

		\begin{itemize}\setl

			\item[1.] $v(P(c_1,\dots,c_n)) = 1$ iff $\langle
			c_1^{\wh{\A}},\dots,c_n^{\wh{\A}} \rangle \in  P_{+}^{\A}$, for every
			$c_1,\dots,c_n \in \mathcal{C}_{\A}$;

			\item[2.] $v(\neg P(c_1,\dots,c_n)) = 1$ iff $\langle
			c_1^{\wh{\A}},\dots,c_n^{\wh{\A}} \rangle \in P_{-}^{\A}$, for every
			$c_1,\dots,c_n \in \mathcal{C}_{\A}$;

			\item[3.] $v(A \land B) = 1$ iff $v(A) = 1$ and $v(B) = 1$;

			\item[4.] $v(A \lor B) = 1$ iff $v(A) = 1$ or $v(B) = 1$;

			\item[5.] $v(\neg (A \land B)) = 1$ iff $v(\neg A) = 1$ or $v(\neg B) = 1$;

			\item[6.] $v(\neg (A \lor B)) = 1$ iff $v(\neg A) = 1$ and $v(\neg B) = 1$;

			\item[7.] $v(\neg \neg A) = 1$ iff $v(A) = 1$;

			\item[8.] If $v(\con A) = 1$, then $v(A) = 1$ iff $v(\neg A) = 0$;

			\item[9.] $v(\con A) = 1$ iff $v(\incon A) = 0$;

			\item[10.] $v(\forall xA) = 1$ iff $v(A(\ov{a}/x)) = 1$, for every $a \in \D$;

			\item[11.] $v(\exists xA) = 1$ iff $v(A(\ov{a}/x)) = 1$, for some $a \in \D$;

			\item[12.] $v(\neg \forall xA) = 1$ iff $v(\neg A(\ov{a}/x)) = 1$, for some $a
			\in \D$;

			\item[13.] $v(\neg \exists xA) = 1$ iff $v(\neg A(\ov{a}/x)) = 1$, for every $a
			\in \D$.

			\item[14.] If $A'$ is an alphabetic variant of $A$, then $v(A') = v(A)$.

			\item[15.] Let $A \in Form(\sig_{\A})$ be such that no variables other than $x$
			are free in $A$, and let $c_{1},c_{2} \in \mathcal{C}_{\A}$. If $c_{1}^{\wh{\A}}
			= c_{2}^{\wh{\A}}$ and $v(A(c_{1}/x)) = v(A(c_{2}/x))$, then $v(\# A(c_{1}/x)) =
			v((\# A(c_{2}/x))$ (where $\# \in  \{\neg, \con, \incon\}$).

			\end{itemize}
		\end{df}

	\begin{df}\label{df:interpretation} Let \sig\ be a signature. An \sig-\textit{interpretation} is
	a pair $\langle \A, v \rangle$ such that \A\ is an \sig-structure and $v$ is an \A-valuation.\\

	\noi A sentence $A$ is said to \textit{hold in} the interpretation $\langle \A, v \rangle$ ($\A,
	v \vDash A$) if and only if $v(A) = 1$; and a set of sentences $\Gamma$ is said to hold in
	$\langle \A, v \rangle$ ($\A, v \vDash \Gamma$) if and only if every element of $\Gamma$ holds in
	$\langle \A, v \rangle$.  $\Gamma$ is said to have a \textit{model} if it holds in some
	interpretation. Finally, $A$ is a \emph{semantic consequence} of $\Gamma$ ($\Gamma \vDash A$) if
	and only if $\A, v \vDash A$ whenever $\A, v \vDash \Gamma$ for every interpretation $\langle \A,
	v \rangle$.
	\end{df}

	\begin{rmrk} \label{rmrk:qletf.inte} \noi Definitions \ref{df:structure}, \ref{df:valuation}, and
	\ref{df:interpretation} above deserve some comments. The definition of \qletf-structures
	resembles very much the corresponding definition in classical first-order logic, except for the
	interpretation given to the predicate letters in terms of extensions and anti-extensions.
	Specifying a \qletf-structure, however, is not sufficient to determine the semantic values of
	\textit{all} sentences. This is a consequence of the fact that in \qletf\ some sentential
	connectives, viz. $\neg$, $\con$,  and $\incon$, are non-deterministic, which means that the
	semantic value of, say, $A$ does not always determine the semantic values of $\neg A$, $\con A$,
	and $\incon A$. For instance, even when $A$ holds in a certain structure \A, there are
	circumstances in which $\con A$ holds in \A, and circumstances in which $\con A$ does not hold in
	\A.

	Regarding Definition \ref{df:valuation}, notice first that valuations assign the values $1$ or
	$0$ to the sentences of the \textit{diagram language} of \sig, which, of course, include all the
	sentences in $Sent(\sig)$. Resorting to diagram languages is required to make sure that the
	quantifiers range over \textit{all} the objects of the domain of \A. Since quantifiers are given
	a  substitutional interpretation -- i.e., the semantic value of a formula $\forall xA$ depends on
	the semantic values of all the substitution instances of $A$ -- we need to extend  the original
	language  with  a new individual constant for each element of the domain $\D$,  and to extend 
	the interpretation function of the
	original structure accordingly. Notice further that clause (14) explicitly requires that any two
	formulas that differ only in some of their bound variables must be assigned the same value by a
	valuation. This clause is the counterpart of rule \textit{AV}, without which formulas such as $\con
	\forall xA$ and $\con \forall yA(y/x)$ (where $y$ does not occur in $A$) cannot be proven to be
	equivalent.\footnote{If void quantifiers were allowed, there would be sentences that are intuitively equivalent, but which could receive different semantic values even in the presence of clause (14) -- e.g., $\con \forall xPc$ and $\con Pc$. That is the reason why 
	we excluded ``formulas" in which void quantifiers occur from the set of formulas. An alternative approach would be to allow for void quantifiers, but extend the definition of alphabetic variants in a such a way that formulas that differ by the presence of one or more void quantifiers would also count as alphabetic variants of one another.}. 
	Finally, clause (15) is required for similar reasons: had it been missing, nothing
	would prevent $\con A(c_{1}/x)$ and $\con A(c_{2}/x)$ from being assigned different values by a
	valuation even though $A(c_{1}/x)$ and $A(c_{2}/x)$ had the same value and $c_{1}$ and $c_{2}$
	were interpreted as the same individual of the domain\footnote{Later on, in Proposition \ref{prop:substitution_01}, we will prove a
	generalization of clause (15).}.
	\end{rmrk}

Let us move on now to the task of proving the completeness of \qletf.

\section{Soundness and Completeness of \qletf}\label{sec:soundnessandcompleteness}

We shall start by establishing that \qletf\ is sound with respect to the class of all \qletf-structures,
leaving the proof of its completeness to the next section.\footnote{The proof of the completeness of
\qletf\ below is based on the one presented in \cite{qmbc} for the logic \qmbc, but requires a few
adjustments in order to comply with specificities of \qletf. In fact, the semantics and the completeness
proof of \qmbc\ in \cite{qmbc} turn out to be special cases of the semantics and the proof presented
here.}

\subsection{Soundness}\label{sec:soundness}

	\begin{prop}\label{prop:substitution_01} Let $\sig = \langle \mathcal{C}, \mathcal{P} \rangle$ be
	a signature and let $ \langle \A, v \rangle$ be an \sig-interpretation. Let $A \in
	Form(\sig_{\A})$ be such that no variables other than $x$ are free in $A$, and $c_{1}, c_{2} \in
	\mathcal{C}_{\A}$. If $c_{1}^{\wh{\A}} = c_{2}^{\wh{\A}}$, then $v(A(c_{1}/x)) = v(A(c_{2}/x))$.
	\end{prop}

	\begin{pr} The result follows by induction on the complexity of $A$ and uses clause (15) of
	Definition \ref{df:valuation}.
	\end{pr}

It is worth noting at this point that if the language of \qletf\ did not include the connectives $\con$
and $\incon$, then the proof of the lemma above would not require clause (15) of Definition
\ref{df:valuation}. In effect, it can proven, \textit{without using that clause}, that:

	\begin{prop}\label{prop:substitution_02} Let $A \in Form(\sig_{A})$ be such that (i) no variables other
	than $x$ are free in $A$ and (ii) $A$ is $\con$- and $\incon$-free. Let $c_{1}, c_{2} \in
	\mathcal{C}_{\A}$.  If $c_{1}^{\wh{\A}} = c_{2}^{\wh{\A}}$, then $v(A(c_{1}/x)) = v(A(c_{2}/x))$.
	\end{prop}

	\begin{pr} The result follows by induction on the complexity of $A$: If $A$ is an atomic formula,
	then the result follows immediately from Definition \ref{df:valuation}(1). If $A$ has the form
	$\neg B$, then there are a few cases : If $B$ is an atomic formula, then the result follows from
	Definition \ref{df:valuation}(2); if $B$ has the form $\neg C$, then the result follows from the
	induction hypothesis (IH) and Definition \ref{df:valuation}(7); if $B$ has the form $A \land B$
	or $A \lor$, then the result follows from (IH) and Definition \ref{df:valuation}(5) and (6); and
	if $B$ has the form $\forall xC$ or $\exists x C$, then the result follows from (IH) and
	Definition \ref{df:valuation}(12) and (13). The remaining cases (viz., $A = B \land C$, $A = B
	\lor C$, $A = \forall xB$, and $A = \exists xB$) are immediate consequences of (IH) and the
	corresponding clauses in Definition \ref{df:valuation} -- i.e., clauses (3), (4), (10), and (11),
	respectively.
	\end{pr}

\noi Moreover, if the language of \qletf\ did not include $\con$ and $\incon$, clause (14) would be
unnecessary as well, for it would then be provable by a straightforward induction on the complexity of
$A$.

The following technical result, which will be used in the proof of the soundness of \qletf, is an
immediate consequence of Proposition \ref{prop:substitution_01} above:

	\begin{cor}\label{cor:substitution_03} Let $\sig = \langle \mathcal{C}, \mathcal{P} \rangle$ be a
	signature such that $c \in \mathcal{C}$. Let $\A = \langle \D, \I \rangle$ be an \sig-structure
	and assume that $a \in \D$. Let $\A' = \langle \D, \I' \rangle$ be the \sig-structure such that
	$\I'$ is just like \I\ except that $\I'(c) = a$. Let $v$ be an \A-valuation and let $v'$ be an
	$\A'$-valuation that agrees with $v$ on all sentences of $Sent(\sig)$ in which $c$ does not
	occur. If no variables other than $x$ are free in $A \in Form(\sig)$ and $c$ does not occur in
	$A$, then $v(A(\ov{a}/x)) = v'(A(c/x)) = 1$.
	\end{cor}

	\begin{pr} Since $\ov{a}^{\wh{\A'}} = c^{\wh{\A'}}$, it follows from Proposition
	\ref{prop:substitution_01} that $v'(A(\ov{a}/x)) = v'(A(c/x))$. But since $c$ does not occur in
	$A$, $v(A(\ov{a}/x)) = v'(A(\ov{a}/x)) = 1$. Therefore, $v(A(\ov{a}/x)) = 1$ if and only if
	$v'(A(c/x)) = 1$.
	\end{pr}

	\begin{tr}\label{tr:soundness} \textbf{(Soundness Theorem)} Let $\sig$ be a signature and $\Gamma
	\cup \{A\} \subseteq Sent(\sig)$. If $\Gamma \vdash A$, then $\Gamma \vDash A$.
	\end{tr}

	\begin{pr} Let $\Theta$ be a derivation of $A$ from $\Gamma$ in \qletf\ and let $n$ be the number
	of nodes in $\Theta$. If $n = 1$, then either $A \in \Gamma$ or $A$ is the result of an
	application of $Comp$ or $\oid I$. If $A \in \Gamma$, then $\Gamma \vDash A$, since $\vDash$ is
	reflexive. If $A = \con B \lor \incon B$ results from an application of $Comp$, then $v(\con B) =
	1$ if and only if $v(\incon B) \neq 1$, for every \sig-interpretation $\langle \A, v \rangle$ (by
	Definition \ref{df:valuation}(8)). Hence, $v(\con B) = 1$ or $v(\incon B) = 1$. By Definition
	\ref{df:valuation}(4), $v(\con B \lor \incon B) = 1$.  Thus, $\A, v \vDash \con B \lor \incon B$,
	and therefore $\vDash A$. If, on the other hand, $A = c \oid c$ results from an application of
	$\oid I$, then, $\langle c^{\wh{\A}},c^{\wh{\A}} \rangle \in \ \oid_{+}^{\A}$, for every
	structure \A\ (by Definition \ref{df:valuation}(5)). By Definition \ref{df:valuation}(1), $v(c
	\oid c) = 1$, for every \A-valuation $v$. Hence, $\A, v \vDash c \oid c$, and therefore $\vDash
	A$.

	\vspace{0.25cm} \noi Suppose that $n > 1$ and that the result holds for every derivation $\Theta'$ with
	fewer nodes than $\Theta$. We shall prove that $\Gamma \vDash A$. Since $n > 1$, $A$ results from an
	application of one of the rules of \qletf\ other than $Comp$ and $\oid I$. In these cases, the result
	follows almost immediately from the lemmas above together with the corresponding clauses of Definition
	\ref{df:valuation}. For instance:

		\begin{itemize}\setl

		\item[3.] Let $A = C \vee \forall xB$ and suppose that it results from an application of
		rule $\forall I$ to $C \vee B(c/x)$. Hence, there is a derivation $\Theta'$ of $C \vee
		B(c/x)$ from $\Gamma$ in \qletf\ such that $\Theta'$ has fewer nodes than $\Theta$. Let
		$\Gamma_{0} \subseteq \Gamma$ be set of hypotheses on which $C \vee B(c/x)$ depends in
		$\Theta'$. Clearly, $\Theta'$ is a derivation of $C \vee B(c/x)$ from $\Gamma_{0}$ and,
		given the restrictions upon $\forall I$, $c$ occurs neither in $B$, $C$ nor in any
		element of $\Gamma_{0}$. By (IH), $\Gamma_{0} \vDash C \vee B(c/x)$. Let $\langle \A, v
		\rangle$ be an arbitrary \sig-interpretation and suppose that $\A, v \vDash \Gamma_{0}$.
		Let $a$ be an arbitrary element of $\D$ and consider then an interpreation $\langle \A',
		v' \rangle$ such that $\A'$ is just like \A\ except that $\I'(c) = a$ and $v'$ agrees
		with $v$ on all sentences in which $c$ does not occur. Since $c$ does not occur in
		$\Gamma_{0}$, $v(D) = 1$ if and only if $v'(D) = 1$, for every $D \in \Gamma_{0}$. Thus,
		$\A', v' \vDash \Gamma_{0}$, and therefore $\A', v' \vdash C \lor B(c/x)$. By Corollary
		\ref{cor:substitution_03}, $v(C(\ov{a}/x) \lor B(\ov{a}/x)) = 1$. Now, since $C(\ov{a}/x)
		= C$, it follows that $v(C \lor B(\ov{a}/x)) = 1$. By Definition \ref{df:valuation}(4),
		$v(C) = 1$ or $v(B(\ov{a}/x)) = 1$. But since $a$ was arbitrary, it follows that $v(C) =
		1$ or $v(B(\ov{a}/x)) = 1$, for every $a \in \D$. By Definition \ref{df:valuation}(10),
		$v(C) = 1$ or $v(\forall xB) = 1$, and by Definition \ref{df:valuation}(4) again, $v(C
		\vee \forall xB) = 1$. That is: $\A, v \vDash A$.

		\item[4.] Let $A = \neg B(c/x)$ and suppose that it results from an application of $\neg
		\exists E$ to $\neg \exists xB$. Hence, there is a derivation $\Theta'$ of $\neg \exists
		xB$ from $\Gamma$ in \qletf\ such that $\Theta'$ has fewer nodes than $\Theta$. By (IH),
		$\Gamma \vDash \neg \exists xB$. Suppose that $\A, v \vDash \Gamma$. Thus, $v(\neg
		\exists xB) = 1$. By Definition \ref{df:valuation}(13), $v(\neg B(\ov{a}/x)) = 1$, for
		every $a \in \D$. Let $b \in \D$ be such that $c^{\wh{\A}} = b$. Therefore, $v(\neg
		B(\ov{b}/x)) = 1$. Since $\ov{b}^{\wh{\A}} = c^{\wh{\A}}$, it follows by Proposition
		\ref{prop:substitution_01} that $v(\neg B(c/x)) = 1$.  Therefore, $\A, v \vDash A$.

		\item[5.] Let $A = B(c_{2}/x)$ and suppose that it results from an application of $\oid
		E$ to $c_{1} \oid c_{2}$ and $B(c_{1}/x)$. Hence, there are derivations $\Theta_{1}$ and
		$\Theta_{2}$ of respectively $c_{1} \oid c_{2}$ and $B(c_{1}/x)$ from $\Gamma$. By (IH),
		$\Gamma \vDash c_{1} \oid c_{2}$ and $\Gamma \vDash B(c_{1}/x)$. Suppose that $\A, v
		\vDash \Gamma$. Thus, $v(c_{1} \oid c_{2}) = v(B(c_1/x)) = 1$. As a result,
		$c_{1}^{\wh{\A}} = c_{2}^{\wh{\A}}$ (by Definition \ref{df:structure}(5)). Since
		$v(B(c_{1}/x)) = 1$, it follows by Proposition \ref{prop:substitution_01} that
		$v(B(c_2/x)) = 1$.

		\end{itemize}

	\vspace{0.25cm} The proof of the remaining cases are left to the reader.
	\end{pr}

\subsection{Completeness}\label{sec:completeness}

Let us now prove the completeness of \qletf. As usual, the proof will be divided up in two separate
steps. First, we shall prove that, given a set of sentences $\Gamma$ that does not prove a certain
sentence $A$, $\Gamma$ can be extended to a set $\Delta$ that (i) still does not prove $A$, (ii) is
closed under $\vdash$, and (iii) has \emph{witnesses} for every universal and existential sentence -- in
the sense that if $B(c/x) \in \Delta$, for every constant $c$, then $\forall xB \in \Delta$, and if
$\exists xB \in \Delta$, then there is an individual constant $c$ such that $B(c/x) \in \Delta$. Second,
we shall prove that, given a set $\Delta$ satisfying (i)-(iii), it is possible to construct a structure
$\A$ and define an $\A$-valuation $v$ such that all and only the sentences belonging to $\Delta$ hold in
$\A$ and $v$.

	\begin{df}\label{df:henkinset} Let $\sig = \langle \mathcal{C}, \mathcal{P} \rangle$ be a
	signature and $\Delta \subseteq Sent(\sig)$. $\Delta$ is a \emph{Henkin set} if and only if for
	every $B \in Sent(\sig)$ and $x \in \mathcal{V}$ (i) $\Delta \vdash \exists xB$ iff $B(c/x)$, for
	some $c \in \mathcal{C}$; and (ii) $\Delta \vdash \forall xB$ iff $B(c/x)$, for every $c \in
	\mathcal{C}$.
	\end{df}

	\begin{df}\label{df:regularset} Let $\sig$ be a signature and $\Delta \cup \{A\} \subseteq
	Sent(\sig)$.  $\Delta$ is a \emph{regular} set if and only if (i) $\Delta$ is
	\textit{non-trivial}: $\Delta \nvdash A$, for some $A \in Sent(\sig)$; (ii) $\Delta$ is
	\textit{closed}: if $\Delta \vdash A$, then $A \in \Delta$, for every $A \in Sent(\sig)$; and
	(iii) $\Delta$ is a \textit{disjunctive set}: if $\Delta \vdash A \lor B$, then $\Delta \vdash A$
	or $\Delta \vdash B$.
	\end{df}

	\begin{lm}\label{lm:dprops} Let $\sig$ be a signature and $\Delta \cup \{A\} \subseteq
	Sent(\sig)$. If $\Delta$ is a regular Henkin set, then:

		\begin{itemize}\setl

		\item[1.] $B \land C \in \Delta$ iff $B \in \Delta$ and $C \in \Delta$;

		\item[2.] $B \lor C \in \Delta$ iff $B \in \Delta$ or $C \in \Delta$;

		\item[3.] $\neg(B \land C) \in \Delta$ iff $\neg B \in \Delta$ or $\neg C \in \Delta$;

		\item[4.] $\neg(B \lor C) \in \Delta$ iff $\neg B \in \Delta$ and $\neg C \in \Delta$;

		\item[5.] $\neg \neg B \in \Delta$ iff $B \in \Delta$;

		\item[6.] $\con B \in \Delta$ iff $\incon B \notin \Delta$;

		\item[7.] If $\con B \in \Delta$, then $B \in \Delta$ iff $\neg B \notin \Delta$;

		\item[8.] If $B'$ is an alphabetic variant of $B$, then $B' \in \Delta$ iff $B \in
		\Delta$;

		\item[9.] $\forall xB \in \Delta$ iff $B(c/x) \in \Delta$, for every $c \in \mathcal{C}$.

		\item[10.] $\exists xB \in \Delta$ iff $B(c/x) \in \Delta$, for some $c \in \mathcal{C}$;

		\item[11.] $\neg \forall xB \in \Delta$ iff $\neg B(c/x) \in \Delta$, for some $c \in
		\mathcal{C}$.

		\item[12.] $\neg \exists xB \in \Delta$ iff $\neg B(c/x) \in \Delta$, for every $c \in
		\mathcal{C}$;

		\end{itemize}
	\end{lm}

	\begin{pr} (1)-(8) are straightforward consequences of the assumption that $\Delta$ is a regular
	set together with the rules of \qletf. (9) and (10) follow immediately from rules $\forall E$ and
	$\exists I$ and the hypothesis that $\Delta$ is a Henkin set. As for (11) and (12), they follow
	respectively from (10) and (9) and the fact that $\neg \forall xB$ and $\exists x \neg B$, and
	$\neg \exists xB$ and $\forall x \neg B$ are derivable from one another in \qletf.
	\end{pr}

We can now prove that if $\Gamma$ is such that $\Gamma \nvdash A$, it can be extended to a regular Henkin
set $\Delta$ such that $\Delta \nvdash A$. For the sake of simplicity, the proof of Lemma
\ref{lm:lindenbaum} (and of Lemma \ref{lm:canonical} and Theorem \ref{tr:completeness}) assumes that
$\sig$ has denumerable non-logical symbols -- and so that the set of formulas generated by $\sig$ is also
denumerable. The reader should, however, encounter no difficulties in generalizing those results to
transfinite languages, and while stating and proving some corollaries of the completeness of \qletf\ at
the end of this section we shall assume that the completeness theorem also holds for such languages.

	\begin{lm}\label{lm:lindenbaum} Let $\mathcal{S} = \langle \mathcal{C}, \mathcal{P} \rangle$ be a
	signature and $\Gamma \cup \{A\} \subseteq Sent(\mathcal{S})$. If $\Gamma \nvdash A$, then there
	is a signature $\sig^{+} = \langle \mathcal{C}^{+}, \mathcal{P} \rangle$ and a regular Henkin set
	$\Delta \subseteq Sent(\sig^{+})$ such that $\mathcal{C} \subseteq \mathcal{C}^{+}$, $\Gamma
	\subseteq \Delta$, and $\Delta \nvdash A$.
	\end{lm}

	\begin{pr} Let $\mathcal{C}^{+} = \mathcal{C} \cup \{c_{i}: i \in \mathbb{N}\}$, with $\{c_{i}: i \in
	\mathbb{N}\} \cap (\mathcal{C} \cup \mathcal{P}) = \emptyset$, and adopt a fixed enumeration
	$B_{0},B_{1},B_{2},\dots$ of the sentences in $Sent(\sig^{+})$. Define the sequence $\langle j_{n}
	\rangle_{n \in \mathbb{N}}$ of natural numbers as follows:

		\begin{itemize}

		\item[1.] $j_{0} =$ the least natural number $k$ such that $c_{k}$ does not occur in
		$B_{0}$.

		\item[2.] $j_{n+1} =$ the least natural number $k$ such that $c_{k}$ does not occur in
		$B_{n+1}$ and for every $i \leq n$, $k \neq j_{i}$.

		\end{itemize}

	\noi Now, consider the sequences $\langle \Gamma_{n} \rangle_{n \in \mathbb{N}}$ and $\langle
	A_{n} \rangle_{n \in \mathbb{N}}$ defined as follows:

		\begin{itemize}

		\item[1.] $\Gamma_{0} = \Gamma$ and $A_{0} = A$;

		\item[2.1.] $\Gamma_{n+1} = \left\{\begin{array}{lll} \Gamma_{n}	& \mbox{if } \Gamma_{n}, B_{n} \vdash A_{n}								\\

											&													\\

			\Gamma_{n} \cup \{B_{n}\} 					& \mbox{if } \Gamma_{n}, B_{n} \nvdash A_{n} \mbox{; and }						\\
											& B_{n} \neq \exists xC \mbox{, for every } x \in \mathcal{V} \mbox{ and  } C \in Form(\sig^{+})	\\

											&													\\

			\Gamma_{n} \cup \{B_{n}, C(c_{j_{n}}/x) \} 			& \mbox{ if } \Gamma_{n}, B_{n} \nvdash A_{n} \mbox{; and }						\\
											& B_{n} = \exists xC \mbox{, for some } x \in \mathcal{V} \mbox{ and } C \in Form(\sig^{+})		\\

											&													\\

		\end{array} \right.$

		\item[2.2.]  $A_{n+1} = \left\{\begin{array}{lll} A_{n}	& \mbox{if } \Gamma_{n}, B_{n} \nvdash A_{n}								\\

									&													\\

			A_{n} \vee B_{n}  				& \mbox{if } \Gamma_{n}, B_{n} \vdash A_{n} \mbox{; and }						\\
									& B_{n} \neq \forall xC \mbox{, for every } x \in \mathcal{V} \mbox{ and  } C \in Form(\sig^{+})	\\

									&													\\

			A_{n} \vee B_{n} \vee C(c_{j_{n}}/x)   		& \mbox{ if } \Gamma_{n}, B_{n} \vdash A_{n} \mbox{; and }						\\
									& B_{n} = \forall xC \mbox{, for some } x \in \mathcal{V} \mbox{ and } C \in Form(\sig^{+})		\\

									&													\\

		\end{array} \right.$

		\end{itemize}

	\vspace{0.25cm}

	\noi Let $\Delta = \bigcup_{n \in \mathbb{N}}\Gamma_{n}$. Clearly, $\Gamma \subseteq \Delta$. We
	shall now prove that $\Delta$ is a regular Henkin set such that $\Delta \nvdash A$. It suffices
	to prove the following facts:

		\begin{itemize}

		\item[1.] For every $n \in \mathbb{N}$, $\Gamma_{n} \nvdash A_{n}$: The proof proceeds by
		induction on $n$.  By the initial hypothesis, $\Gamma_{0} \nvdash A_{0}$. Suppose that
		$\Gamma_{n} \nvdash A_{n}$ (IH). There are two cases: either (I) $\Gamma_{n}, B_{n}
		\vdash A_{n}$ or (II) $\Gamma_{n}, B_{n} \nvdash A_{n}$. If (I) $\Gamma_{n}, B_{n} \vdash
		A_{n}$, then $\Gamma_{n+1} = \Gamma_{n}$. Suppose then that $\Gamma_{n+1} \vdash
		A_{n+1}$. If $B_{n} \neq \forall xC$, then $A_{n+1} = A_{n} \vee B_{n}$. Hence,
		$\Gamma_{n} \vdash A_{n} \vee B_{n}$. But since $\Gamma_{n}, A_{n} \vdash A_{n}$ and
		$\Gamma_{n}, B_{n} \vdash A_{n}$ (by (I)), $\Gamma_{n} \vdash A_{n}$, which contradicts
		(IH).  If, on the other hand, $B_{n} = \forall xC$, then $A_{n+1} = A_{n} \vee B_{n} \vee
		C(c_{j_{n}}/x)$, and so $\Gamma_{n} \vdash A_{n} \vee B_{n} \vee C(c_{j_{n}}/x)$. Since
		$c_{j_{n}}$ does not occur in $\Gamma_{n} \cup \{A_{n}, B_{n}, C\}$, we can apply  rule
		$\forall I$ to obtain $\Gamma_{n} \vdash A_{n} \vee B_{n}$:

			\begin{center}

				\bottomAlignProof
								\AxiomC{$(A_{n} \vee B_{n}) \vee C(c_{j_{n}}/x)$}
							\RightLabel{$\forall I$}
							\UnaryInfC{$(A_{n} \vee B_{n}) \vee \forall xC$}
								\AxiomC{$[A_{n} \vee B_{n}]_{1}$}
								\AxiomC{$[B_{n}]_{1}$}
							\RightLabel{$\vee I$}
							\UnaryInfC{$A_{n} \vee B_{n}$}
					\RightLabel{$\vee E_{1}$}
					\TrinaryInfC{$A_{n} \vee B_{n}$}
				\DisplayProof

			\end{center}

		\noindent But this result contradicts (IH).

		Now, if (II) $\Gamma_{n}, B_{n} \nvdash A_{n}$, then $A_{n+1} = A_{n}$ and either
		$\Gamma_{n+1} = \Gamma \cup \{B_{n}\}$ or $\Gamma_{n+1} = \Gamma_{n} \cup \{\exists xC,
		C(c_{j_{n}}/x)\}$. Suppose one more time that $\Gamma_{n+1} \vdash A_{n+1}$. If
		$\Gamma_{n+1} = \Gamma_{n} \cup \{B_{n}\}$, then $\Gamma_{n}, B_{n} \vdash A_{n}$, which
		contradicts (II). Hence, $\Gamma_{n+1} = \Gamma_{n} \cup \{\exists xC, C(c_{j_{n}}/x)\}$.
		Since $c_{j_{n}}$ does not occur in $\Gamma_{n} \cup \{A_{n}, B_{n}, C\}$, we can apply
		rule $\exists E$ to obtain $\Gamma_{n}, B_{n} \vdash A_{n}$:

		\begin{center}

			\bottomAlignProof
						\AxiomC{$\exists xC (= B_{n})$}
						\AxiomC{$[C(c_{j_{n}}/x)]_{1}$} \noLine
					\UnaryInfC{$\vdots$} \noLine
					\UnaryInfC{$A_{n}$}
				\RightLabel{$\exists E_{1}$}
				\BinaryInfC{$A_{n}$}
			\DisplayProof

		\end{center}

		\noindent But this result contradicts (II).

		\item[2.] For every $n \in \mathbb{N}$, $\Delta \nvdash A_{n}$ (in particular, $\Delta
		\nvdash A$): Suppose that $\Delta \vdash A_{n}$, for some $n \in \mathbb{N}$. Hence,
		$\Gamma_{m} \vdash A_{n}$, for some $m \in \mathbb{N}$. If $m \leq n$, then $\Gamma_{n}
		\vdash A_{n}$, since $\Gamma_{m} \subseteq \Gamma_{n}$; and if $m > n$, then $\Gamma_{m}
		\vdash A_{m}$, given that $A_{m} = A_{n} \vee C$, for some $C \in Sent(\sig^{+})$. In
		either case, there is a  contradiction with item (1) above.

		\item[3.] If $\Delta \vdash C$, then $C \in \Delta$: Suppose that $\Delta \vdash C$ and
		that $C \notin \Delta$. Let $n \in \mathbb{N}$ be such that $C = B_{n}$. Since $B_{n}
		\notin \Delta$, $B_{n} \notin \Gamma_{n+1}$. Hence, $\Gamma_{n}, B_{n} \vdash A_{n}$, and
		so $\Delta, B_{n} \vdash A_{n}$.  Therefore, $\Delta \vdash A_{n}$, which contradicts (2)
		above.

		\item[4.] If $\Delta \vdash C \vee D$,  then $\Delta \vdash C$ or $\Delta \vdash D$:
		Suppose that $\Delta \vdash C \vee D$ and that $\Delta \nvdash C$ and $\Delta \nvdash D$.
		Thus, $C \notin \Delta$ and $D \notin \Delta$. Let $m, n \in \mathbb{N}$ be such that
		$B_{m} = C$ and $B_{n} = D$.  Hence, $B_{m} \notin \Gamma_{m+1}$ and $B_{n} \notin
		\Gamma_{n+1}$, and so $\Gamma_{m}, B_{m} \vdash A_{m}$ and $\Gamma_{n}, B_{n} \vdash
		A_{n}$. Since $\Gamma_{m}, \Gamma_{n} \subseteq \Delta$, it follows that $\Delta, B_{m}
		\vee B_{n} \vdash A_{m} \vee A_{n}$. Now, if $m \leq n$, then $A_{m} \vee A_{n} \vdash
		A_{n}$ (by the way $\langle A_{n} \rangle_{n \in \mathbb{N}}$ was defined).  Hence,
		$\Delta, B_{m} \vee B_{n} \vdash A_{n}$. But since $\Delta \vdash B_{m} \vee B_{n}$,
		$\Delta \vdash A_{n}$, which contradicts (2) above. If $n < m$, then $A_{m} \vee A_{n}
		\vdash A_{m}$, and the result follows in an analogous way.

		\item[5.] For every $C \in Sent(\sig^{+})$, $\Delta \vdash \forall xC$ if and only if
		$\Delta \vdash C(c/x)$, for every $c \in \mathcal{C}^{+}$: It suffices to prove that if
		$\Delta \vdash C(c/x)$, for every $c \in \mathcal{C}^{+}$, then $\Delta \vdash \forall
		xC$, since the other direction is an immediate consequence of rule $\forall E$. We shall
		prove the contrapositive: Let $n \in \mathbb{N}$ be such that $\forall xC = B_{n}$ and
		suppose that $\Delta \nvdash B_{n}$.  Hence, $\Gamma_{n+1} \nvdash B_{n}$, and so $B_{n}
		\notin \Gamma_{n+1}$. By the definition of $\langle \Gamma_{n} \rangle_{n \in
		\mathbb{N}}$, it then follows that $\Gamma_{n}, B_{n} \vdash A_{n}$ and that $A_{n+1} =
		A_{n} \vee B_{n} \vee C(c_{j_{n}}/x)$. Suppose that $\Delta \vdash C(c_{j_{n}}/x)$.
		Thus, $\Delta \vdash A_{n+1}$, which contradicts (2) above. Hence, $\Delta \nvdash
		C(c/x)$ for at least one $c \in \mathcal{C}^{+}$.

		\item[6.] For every $C \in Sent(\sig^{+})$, $\Delta \vdash \exists xC$ if only if $\Delta
		\vdash C(c/x)$, for some $c \in \mathcal{C}^{+}$: As in (5), we shall only prove the
		(contrapositive of the) left-to-right direction, since the other direction follows
		immediately from rule $\exists I$.  Let $n \in \mathbb{N}$ be such that $\exists xC =
		B_{n}$ and suppose that $\Delta \nvdash C(c/x)$, for every $c \in \mathcal{C}^{+}$. In
		particular, $\Delta \nvdash C(c_{j_{n}}/x)$. Hence, $\Gamma_{n+1} \nvdash C(c_{j_{n}}/x)$
		and, by construction, $\Gamma_{n}, B_{n} \vdash A_{n}$. By the monotonicity of $\vdash$,
		it then follows that $\Delta, B_{n} \vdash A_{n}$. Suppose that $\Delta \vdash B_{n}$. By
		the transitivity of $\vdash$, $\Delta \vdash A_{n}$, which contradicts (2) above. Hence,
		$\Delta \nvdash  B_{n}$ (i.e., $\Delta \nvdash \exists xC$).

		\end{itemize}
	\end{pr}

	\begin{lm}\label{lm:identity} Let $\sig$ be a first-order structure, $\Gamma \subseteq
	Sent(\sig)$, and $c_{1},c_{2},c_{3} \in \mathcal{C}$. Then:

		\begin{itemize}\setl

		\item[1.] $\Gamma \vdash c_{1} \oid c_{1}$;

		\item[2.] If $\Gamma \vdash c_{1} \oid c_{2}$, then $\Gamma \vdash c_{2} \oid c_{1}$; and

		\item[3.] If $\Gamma \vdash c_{1} \oid c_{2}$ and $\Gamma \vdash c_{2} \oid c_{3}$, then
		$\Gamma \vdash c_{1} \oid c_{3}$.

		\end{itemize}
	\end{lm}

	\begin{pr} (1)-(3) result immediately from applying rules $\oid I$ and $\oid E$.
	\end{pr}

The only remaining step to finish the proof of the completeness of \qletf\ is to show that, given a
regular Henkin set $\Delta$ that does not prove $A$, one can construct a \emph{canonical model} $\A$ and
a valuation $v$ in $\A$ such that all elements of $\Delta$ (and no others) hold in $\A$ and $v$. Since $A
\notin \Delta$, this result will then be enough to conclude that $\Delta \nvDash A$ -- which, when
combined with Lemma \ref{lm:lindenbaum}, suffices for proving the completeness of \qletf.

Lemma \ref{lm:canonical} below has a rather long proof and so it might be worthy detailling its
structure. The proof comprises three different parts: in the first part we define a canonical structure
$\A$ that interprets the non-logical symbols of the relevant language. As usual, the domain of $\A$ will
be composed by the terms of the language (or rather, by certain equivalence classes thereof), while the
non-logical symbols will be interpreted in terms of the \emph{derivability-from-}$\Delta$ relation. As a
result, for each constant $c \in \mathcal{C}$, the corresponding diagram language will include a new
individual constant $\ov{[c]}$, for each constant of the original language. In the second part we define
a valuation $v$ such that all and only the elements of $\Delta$ hold in $\A$ and $v$. Finally, in the
third part, we prove that $v$, thus defined, is indeed a valuation in $\A$, which comes down to showing
that it satisfies all clauses of Definition \ref{df:valuation}.

	\begin{lm}\label{lm:canonical} Let $\sig = \langle \mathcal{C}, \mathcal{P} \rangle$ be a
	signature and $\Delta \cup \{A\} \subseteq Sent(\sig)$. If $\Delta$ is a regular Henkin set, then
	there is an \sig-interpretation $\langle \A, v \rangle$ such that $\A, v \vDash \Delta$ and $\A,
	v \nmodels A$.
	\end{lm}

	\begin{pr} Let $\sim$ be the relation on $\mathcal{C}$ defined by: $c_{1} \sim c_{2}$ iff $\Delta
	\vdash c_{1} \oid c_{2}$. For each $c \in \mathcal{C}$, let $[c] = \{c' \in \mathcal{C}: c' \sim
	c\}$. Notice that since $\sim$ is an equivalence relation (by Lemma \ref{lm:identity}), $[c_{1}]
	= [c_{2}]$ if, and only if, $c_{1} \sim c_{2}$. Now, define the \sig-structure $\A = \langle \D,
	\I \rangle$ as follows:

		\begin{itemize}\setl

		\item[1.] $\D = \{[c]: c \in \mathcal{C}\}$\footnote{$\D$ is thus the partition of
		$\mathcal{C}$ generated by $\sim$.};

		\item[2.] For every $c \in \mathcal{C}$, $c^{\A} = [c]$;

		\item[3.] For every $n$-ary predicate letter $P \in \mathcal{P}$, and $c_{1},\dots,c_{n}
		\in \mathcal{C}$:

			\begin{itemize}\setl

			\item[3.1.] $\langle [c_{1}],\dots,[c_{n}] \rangle \in P_{+}^{\A}$ iff $\Delta
			\vdash P(c_{1},\dots,c_{n})$; and

			\item[3.2.] $\langle [c_{1}],\dots,[c_{n}] \rangle \in P_{-}^{\A}$ iff $\Delta
			\vdash \neg P(c_{1},\dots,c_{n})$.

			\end{itemize}

		\end{itemize}

	\noi The following facts guarantee that $\A$ is indeed an \sig-structure:

		\begin{itemize}\setl

		\item[1.] For every $P \in \mathcal{P}_{n}$, the definition of $P^{\A} = \langle
		P^{\A}_{+}, P^{\A}_{-} \rangle$ does not depend on the representatives
		$c_{1},\dots,c_{n}$.

		\item[2.] $\langle [c], [c] \rangle \in \mathcal{\oid}_{+}^{\A}$, for every $c \in
		\mathcal{C}$.

		\end{itemize}

	\noi In order to prove (1), it suffices to show that $\Delta \vdash P(c_{1},\dots,c_{n})$ if, and
	only if, $\Delta \vdash P(c_{1}',\dots,c_{n}')$, whenever $c_{i} \sim c_{i}'$, for every $1 \leq
	i \leq n$ (and similarly for $\neg P(c_{1},\dots,c_{n})$). But this is an immediate consequence
	of the application of rule $\oid E$. (2) follows immediately by rule $\oid I$, clause (3) in the
	definition of $\A$, and the definition of $\sim$.

	Now, define the mapping $*: Term(\sig_{\A}) \longrightarrow Term(\sig)$ as follows: if $t =
	\ov{[c']}$, for some $c' \in \mathcal{C}$, then $c^{*} = c'$; and if $t \in \mathcal{V} \cup
	\mathcal{C}$, then $t^{*} = t$. The mapping $*$ can then be naturally extended to the formulas of
	$\sig_{\A}$ as follows:

		\begin{itemize}\setl

		\item[1.] If $B = P(t_{1},\dots,t_{n})$, then $B^{*} = P(t_{1}^{*},\dots,t_{n}^{*})$;

		\item[2.] If $B = \#C$ ($\# \in \{\neg, \con, \incon\}$), then $B^{*} = \#C^{*}$;

		\item[3.] If $B = C \# D$ ($\# \in \{\land, \lor\}$), then $B^{*} = C^{*} \# D^{*}$;

		\item[4.] If $B = QxC$ ($Q \in \{\forall, \exists\}$), then $B^{*} = QxC^{*}$.

		\end{itemize}

	\noi Notice that ($\dagger$) for every constant $c \in \mathcal{C}_\A$, $c^{\wh{\A}} = [c^{*}]$.
	For if $c = \ov{[c']}$, for some $c' \in \mathcal{C}$, then $\ov{[c']}^{\wh{\A}} = [c'] =
	[c^{*}]$; and if $c \in \mathcal{C}$, then $c^{\wh{\A}} = [c] = [c^{*}]$. Notice further that
	($\dagger \dagger$) $B(c/x)^{*} = B^{*}(c^{*}/x)$, for every $B \in Form(\sig_{\A})$, which can
	be proven by induction on the complexity of $B$.

	Finally, define the mapping $v: Sent(\sig_{\A}) \longrightarrow \{0,1\}$ by: $v(B) = 1$ iff
	$B^{*} \in \Delta$. Given $B \in Sent(\sig_{\A})$, $v$ assigns to $B$ the value $1$ if the
	corresponding formula $B^{*}$ of $Sent(\sig)$ belongs to $\Delta$ (and $0$
	otherwise)\footnote{The mapping $*$ is necessary to ensure that the sentences of the diagram
	language of \sig\ which do not belong to $Sent(\sig)$ get assinged one of the values 1 or 0 by
	$v$. Since the domain of $v$ is $Sent(\sig_{\A})$, we are prevented from defining it as the
	characteristic function of $\Delta$ (all of whose elements belong to $Sent(\sig)$). By making use
	of $*$, we are nonetheless able to define $v$ in terms of \textit{membership-in-}$\Delta$, for
	the values of the elements of $Sent(\sig_{\A})/Sent(\sig)$ are then determined by whether their
	$*$-translations belong to $\Delta$. So, for example, the value of $P(\ov{[c]})$, which is a
	sentence of the diagram language of \sig, is $1$ if, and only if, its $*$-translantion $Pc$
	belongs to $\Delta$.}. Because $A \notin \Delta$ (by hypothesis), $\A$ and $v$ are such that $\A,
	v \vDash \Delta$ and $\A, v \nmodels A$, as required.

	We haven't finished yet, though. For we are still required to show that $v$, as defined above, is
	indeed an \A-valuation -- given that nothing so far guarantees that it satisfies all clauses of
	Definition \ref{df:valuation}. In the remainder of this proof, we shall prove that $v$ satisfies
	as least some of those clauses, leaving the remaining cases to the reader.

		\begin{itemize}\setl

		\item[1.] Let $B \in Sent(\sig_\A)$ be the formula $P(c_{1},\dots,c_{n})$. Then:

			\begin{center}
				\seta
				{\small \begin{tabular}{llll}
				$v(B) = 1$	& iff	& $P(c_{1},\dots,c_{n})^{*} \in \Delta$						&		\\
						& iff	& $P(c_{1}^{*},\dots,c_{n}^{*}) \in \Delta$					&		\\
						& iff	& $\Delta \vdash P(c_{1}^{*},\dots,c_{n}^{*})$					& 		\\
						& iff	& $\langle [c_{1}^{*}],\dots,[c_{n}^{*}] \rangle \in P_{+}^{\A}$		&		\\
						& iff	& $\langle c_{1}^{\wh{\A}},\dots,c_{n}^{\wh{\A}} \rangle \in P_{+}^{\A}$	& ($\dagger$) above	\\
				\end{tabular}}
			\end{center}

		\item[3.] Let $B = (C \land D)$. Then:

			\begin{center}
				\seta
				{\small \begin{tabular}{llll}
				$v(B) = 1$	& iff	& $(C \land D)^{*} \in \Delta$				&				\\
						& iff	& $C^{*} \in \Delta \mbox{ and } D^{*} \in \Delta$	& Lemma \ref{lm:dprops}(1)	\\
						& iff	& $v(C) = 1 \mbox{ and } v(D) = 1$			& 				\\
				\end{tabular}}
			\end{center}

		\item[9.] Suppose that $B = \con C$ and that $v(B) = 1$. Hence, $(\con C)^{*} \in
		\Delta$, and so $\con C^{*} \in \Delta$. By Lemma \ref{lm:dprops}(7), $C^{*} \in \Delta$
		if, and only if, $\neg C^{*} \notin \Delta$. Thus, $v(C) = 1$ if, and only if, $v(\neg C)
		= 0$.

		\item[10.] Let $B = \forall xC$. Then:

			\begin{center}\seta
				{\small \begin{tabular}{llll}

				$v(B) = 1$	& iff	& $(\forall xC)^{*} \in \Delta$							&				\\
						& iff	& $\forall xC^{*} \in \Delta$							& 				\\
					     	& iff	& $C^{*}(c/x) \in \Delta$, for every $c \in \mathcal{C}$			& Lemma \ref{lm:dprops}(9)	\\
						& iff	& $C^{*}(\ov{c^{\wh{\A}}}^{*}/x) \in \Delta$, for ever $c \in \mathcal{C}$ 	& $\ov{c^{\wh{\A}}}^{*} = c$, for every $c \in \mathcal{C}$	\\
						& iff	& $C(\ov{c^{\wh{\A}}}/x)^{*} \in \Delta$, for every $c \in \mathcal{C}$	& ($\dagger \dagger$)	 			\\
					     	& iff	& $C(\ov{a}/x)^{*} \in \Delta$, for every $a \in \D$  			& 				\\
						& iff	& $v(C(\ov{a}/x)) = 1$, for every $a \in \D$					&				\\
				\end{tabular}}
			\end{center}

		\item[11.] Let $B = \exists xC$ and let $c \in \mathcal{C}$ be such that $\exists xC^{*}
		\in \Delta$ if, and only if, $C^{*}(c/x) \in \Delta$, which we know to exist due to Lemma
		\ref{lm:dprops}(10). We then have:

			\begin{center}\seta

				{\small \begin{tabular}{llll}
				$v(B) = 1$	& iff	& $(\exists xC)^{*} \in \Delta$				&				\\
						& iff	& $\exists xC^{*}$					&				\\
						& iff	& $C^{*}(c/x) \in \Delta$				&				\\
						& iff	& $C^{*}(\ov{c^{\wh{\A}}}^{*}/x) \in \Delta$		& $\ov{c^{\wh{\A}}}^{*} = c$	\\
						& iff	& $C(\ov{c^{\wh{\A}}}/x)^{*} \in \Delta$		& ($\dagger \dagger$) 				\\
						& iff	& $C(\ov{a}/x)^{*} \in \Delta$, for some $a \in \D$ 	& 				\\
						& iff	& $v(C(\ov{a}/x)) = 1$, for some $a \in \D$		&				\\
				\end{tabular}}

			\end{center}

		\item[15.] Let $B \in Form(\sig_{\A})$ be such that $x$ is the only variable free in $B$,
		and let $c_1,c_2 \in \mathcal{C}_\A$. Suppose that $c_1^{\wh{\A}} = c_2^{\wh{\A}}$ and
		that $v(B(c_{1}/x)) = v(B(c_2/x))$. Since $c_1^{\wh{\A}} = c_2^{\wh{\A}}$, it follows by
		($\dagger$) that $[c_1^*] = [c_2^*]$. Hence, $\Delta \vdash c_1^{*} \oid c_2^{*}$. Let
		$\# \in \{\neg, \con, \incon\}$. Thus:

			\begin{center}\seta

				{\small \begin{tabular}{llll}

				$(\# B(c_{1}/x))^{*} \in \Delta$	& iff	& $\#(B(c_{1}/x)^{*}) \in \Delta$	&					\\
									& iff	& $\#B^{*}(c_{1}^{*}/x) \in \Delta$	& ($\dagger \dagger$)					\\
									& iff	& $\# B^{*}(c_2^*/x) \in \Delta$	& $\Delta \vdash c_1^* \oid c_2^*$	\\
									& iff	& $\#(B(c_{2}/x)^{*}) \in \Delta$	& ($\dagger \dagger$)					\\
									& iff	& $(\# B(c_2/x))^* \in \Delta$		&					\\
				\end{tabular}}

			\end{center}

		\noi Therefore, $v(\# B(c_{1}/x)) = v(\# B(c_2/x))$.

		\end{itemize}
	\end{pr}

Having proved lemmas \ref{lm:lindenbaum} and \ref{lm:canonical}, proving the completeness of \qletf\ is
pretty much straightforward:

	\begin{tr}\label{tr:completeness} \textbf{(Completeness Theorem)} Let $\sig = \langle
	\mathcal{C}, \mathcal{P} \rangle$ be a signature and $\Gamma \cup \{A\} \subseteq Sent(\sig)$. If
	$\Gamma \vDash A$, then $\Gamma \vdash A$.
	\end{tr}

	\begin{pr} Suppose that $\Gamma \nvdash A$. By Lemma \ref{lm:lindenbaum}, there is a signature
	$\sig^{+} = \langle \mathcal{C}^{+}, \mathcal{P} \rangle$ and a set $\Delta \subseteq
	Sent(\sig^{+})$ such that $\mathcal{C} \subseteq \mathcal{C}^{+}$, $\Gamma \subseteq \Delta$, and
	$\Delta$ is a regular Henkin set that does not prove $A$. By Lemma \ref{lm:canonical}, there
	exists an $\sig^{+}$-interpretation $\langle \A, v \rangle$ such that $\A, v \vDash \Delta$ and
	$\A, v \nmodels A$. Let $\A_{0}$ be the $\sig$-reduct of $\A$ and $v_{0}$ be the restriction of
	$v$ to $Sent(\sig_{\A})$. Clearly, $\A_{0}, v_{0} \vDash B$ if, and only if, $\A, v \vDash B$,
	for every $B \in Sent(\sig)$. As a result, $\A_{0}, v_{0} \vDash \Gamma$ (since $\A, v \vDash
	\Gamma \subseteq \Delta$) and $\A_{0}, v_{0} \nmodels A$. Therefore, $\Gamma \nmodels A$.
	\end{pr}

Compactness and Lowenheim-Skolem theorems are immediate consequences of the soundness and the
completeness of \qletf:

	\begin{cor}\label{cor:compactness} \textbf{(Compactness Theorem)} Let $\Gamma \cup \{A\}
	\subseteq Sent(\sig)$. Then:

		\begin{itemize}\setl

		\item[1.] $\Gamma \vDash A$ if and only if there is a finite subset $\Gamma_{0}$ of
		$\Gamma$ such that $\Gamma_{0} \vDash A$;

		\item[2.] $\Gamma$ has a model if and only if every finite subset $\Gamma_{0}$ of
		$\Gamma$ has a model.

		\end{itemize}
	\end{cor}

	\begin{pr} (1) is an immediate consequence of the soundness and completeness theorems and the
	fact that derivations are finite. As for (2), we shall only prove the right-to-left direction.
	Suppose that every finite subset $\Gamma_{0}$ of $\Gamma$ has a model and that $\Gamma$ does not
	have a model. Hence, $\Gamma \vDash A$, for every sentence $A$, and so there is a sentence $B$
	such that $\Gamma \vDash \con B \land B \land \neg B$. By Corollary \ref{cor:compactness}, it
	then follows that some finite $\Gamma_{0} \subseteq \Gamma$ is such that $\Gamma_{0} \vDash \con
	B \land B \land \neg B$. Therefore, $\Gamma_{0}$ is trivial and does have a model, which
	contradicts the initial hypothesis.
	\end{pr}

	\begin{cor}\label{cor:downlowenheim-skolem} \textbf{(Downward Lowenh\"eim-Skolem Theorem)} Let
	$\sig = \langle \mathcal{C}, \mathcal{P} \rangle$ be a signature whose cardinality is $\lambda$,
	and suppose that $\Gamma \subseteq Sent(\sig)$. If $\Gamma$ has a model, then $\Gamma$ has a
	model whose cardianlity is less than or equal to $\lambda$.

	\end{cor}

	\begin{pr} Since $\Gamma$ has a model, it is non-trivial (by soundness). Hence, there is a
	sentence $A$ such that $\Gamma \nvdash A$. By Lemma \ref{lm:lindenbaum}, there is a signature
	$\sig^{+} = \langle \mathcal{C}^{+}, \mathcal{P} \rangle$, with $|\mathcal{C}^{+}| = \lambda$,
	and a set $\Delta \subseteq Sent(\sig^{+})$ such that $\Delta$ is a regular Henkin set, $\Gamma
	\subseteq \Delta$, and $\Delta \nvdash A$. By Lemma \ref{lm:canonical}, there is an
	$\sig^{+}$-structure $\A^{+}$ whose domain $\D^{+}$ is the set $\{[c]: c \in \mathcal{C}^{+}\}$,
	and there is a valuation $v^{+}$ in $\A^{+}$ such that $\A^{+}, v^{+} \vDash \Delta$.  Since
	$\D^{+}$ is a partition of $\mathcal{C}^{+}$, $|\D^{+}| \leq |\mathcal{C}^{+}|$, and so $|\D^{+}|
	\leq \lambda$.  Finally, let \A\ be the \sig-reduct of $\A^{+}$ and $v$ be the restriction of
	$v^{+}$ to $Sent(\sig_{\A})$. Since $\Gamma \subseteq \Delta$ and the domain $\D$ of \A\ is equal
	to $\D^{+}$, it follows that $\A, v \vDash \Gamma$ and $|\D| \leq \lambda$.
	\end{pr}

	\begin{cor}\label{cor:uplownheim-skolem} \textbf{(Upward Lowhenh\"eim-Skolem Theorem)} Let $\sig
	= \langle \mathcal{C}, \mathcal{P} \rangle$ be a signature whose cardinality is $\lambda$, and
	let $\Gamma \subseteq Sent(\sig)$ be such that $\Gamma \vDash \forall x \forall y \con(x \oid
	y)$. If $\Gamma$ has an infinite model, then $\Gamma$ has a model of cardinality $\kappa$, for
	every $\kappa \geq \lambda$.
	\end{cor}

	\begin{pr} Let $\kappa \geq \lambda$ and let $\sig^{+} = \langle \mathcal{C}^{+}, \mathcal{P}
	\rangle$ be such that $\mathcal{C}^{+} = \mathcal{C} \cup \{c_{\alpha}: \alpha < \kappa\}$.
	Consider the set: $$\Delta = \Gamma \cup \{c_{\alpha} \noid c_{\beta}: \alpha, \beta < \kappa
	\mbox{ and } \alpha \neq \beta\}$$

	\noi Since $\Gamma$ has a model, so does $\Delta$. For let $\langle \A, v \rangle$ be an infinite
	\sig-interpretation such that $\A, v \vDash \Gamma$, and consider an arbitrary finite subset
	$\Delta_{0}$ of $\Delta$. Define $\A^{+} = \langle \D, \I^{+} \rangle$ to be the extension of \A\
	such that:

		\begin{itemize}\setl

		\item[1.] For every $\alpha, \beta < \kappa$ such that $c_{\alpha}$ and $c_{\beta}$ occur in
		$\Delta_{0}$, if $\alpha \neq \beta$, then $\I^{+}(c_{\alpha}) \neq \I^{+}(c_{\beta})$;

		\item[2.] For every $\alpha$ such that $c_{\alpha}$ does not occur in $\Delta_{0}$,
		$\I^{+}(c_{\alpha})$ is a fixed element $a$ of $\D$;

		\item[3.] $I^{+}(\oid) = I(\oid)$.

		\end{itemize}

	\noi Now, extend $v$ to a valuation $v^{+}$ such that for every $c_{1},\dots,c_{n} \in
	\mathcal{C}$ and $c'_{1},\dots,c'_{n} \in \mathcal{C}^{+}$, and for every sentence $A \in
	Sent(\sig)$, if $\I^{+}(c_{i}) = \I^{+}(c'_{i})$, then $v^{+}(A(c'_{1}/c_{1};\dots;c'_{n}/c_{n})
	= v(A)$.  Clearly, $\Delta_{0}$ holds in $\langle \A^{+}, v^{+} \rangle$, since $\langle a, b
	\rangle \in \ \oid^{\A^{+}}_{-}$ if and only if $a \neq b$. Therefore, every finite subset of
	$\Delta$ has a model, and, by corollaries \ref{cor:compactness}(2) and
	\ref{cor:downlowenheim-skolem}, $\Delta$ has a model whose cardinality $\kappa'$ is less than or
	equal to $\kappa$\footnote{Notice that the cardinality of $\sig^{+}$ is $\kappa$ rather than
	$\lambda$.}. But because $\Delta$ includes every sentence $c_{\alpha} \noid c_{\beta}$ and
	$\Delta \vDash \forall x \forall y\con (x \oid y)$, $\kappa'$ must be equal to $\kappa$ (for if
	$\A, v \vDash \con(c_{1} \oid c_{2}) \land c_{1} \noid c_{2}$, then $\A, v \nvDash c_{1} \oid
	c_{2}$, and so $\I(c_{1}) \neq \I(c_{2})$).
	\end{pr}

	\begin{rmrk}\label{rmrk:openformulas} \textit{On extended valuations}
	
	 \noi So far we have defined
	all relevant syntactic and semantic notions with respect to sentences, completely disregarding
	open formulas.  In particular, while presenting the natural deduction system for \qletf\ we have
	replaced the more traditional quantifier rules by corresponding rules in which constants play the
	roles of variables or terms. This choice led us to assume that every language has an infinite
	stock of individual constants to ensure that there will always be enough constants to meet the
	restrictions upon some of the quantifier rules (viz., $\forall I$, $\exists E$, $\neg \forall E$,
	and $\neg \exists I$).

	Now, although focusing on sentences brings some significant technical simplifications, we could
	have adopted a more traditional approach, formulating the natural deduction system for \qletf\
	with the usual rules, and defining the semantic consequence relation to include both open and
	closed formulas. This could be achieved by making use of \textit{extended valuations}
	\cite[see][Def. 7.3.10]{carco.book}, wich can be defined as follows: Given an \sig-interpretation
	$\langle \A, v \rangle$, the extension of $v$ in \A\ is the mapping $\ov{v}: Form(\sig_{\A})
	\times \D^{\mathcal{V}} \longrightarrow \{0,1\}$ such that: $$\ov{v}(A, s) =
	v(A(\ov{s(x_{1})},\dots,\ov{s(x_{n})}/x_{1},\dots,x_{n}))$$

	\noi where $s$ is an assignment of elements of $\D$ to the individual variables and all variables
	free in $A \in Form(\sig)$ are among $x_{1},\dots,x_{n}$.

	Hence, the extension $\ov{v}$ of $v$ assings to an open formula $A(x_{1},\dots,x_{n})$ the value
	assined by $v$ to the sentence $A(\ov{s(x_{1})},\dots,\ov{s(x_{n})})$ -- which results from $A$
	by replacing the variables $x_{1},\dots,x_{n}$ by the constants
	$\ov{s(x_{1})},\dots,\ov{s(x_{n})}$ of the corresponding diagram language.

	This strategy allows us to mimic the definitions usually found in traditional formulations of a
	referential semantics for, say, first-order classical logic. In particular, it can be proven that
	each clause of Definition \ref{df:valuation} can be rewritten in terms of extended valuations.
	For instance, clause (4), together with definition of $\ov{v}$, allows us to prove:

		$$\ov{v}(B \lor C, s) = 1 \mbox{ iff } \ov{v}(B, s) = 1 \mbox{ or } \ov{v}(C, s) = 1$$

	\noi while clause (10) allows us to prove:

		$$\ov{v}(\forall xB, s) = 1 \mbox{ iff } \ov{v}(B, s_{x}^{a}) = 1, \mbox{ for every } a
		\in \D$$

	\noi (where $s_{x}^{a}$ is the assignment that differs from $s$ at most by assigning $a$ to $x$).

	Had we chosen to adopt this strategy and defined the semantic consequence relation accordingly,
	we would also be capable of proving all the results above -- though the corresponding definitions
	and proofs would become much more cumbersome. That this can be done suffices to ensure that
	nothing in this paper hinges on our choice to focus entirely on sentences, and that we could have
	done without the assumption that languages must always have infinitely many individual constants.

	\end{rmrk}

	\begin{rmrk}\textit{On the non-classical identity of \qletf} \label{rmk:identity}

	\noi We remarked in the introduction (p.~\pageref{page.info}) that \qletf\ can be interpreted in
	terms of information, which may be positive or negative, reliable or unreliable. Now, let us
	illustrate the semantics of identity based on this interpretation. 
	
	 There are four scenarios of
	absence of reliable information (i.e. when $v(\cons (a\oid b)) = 0$),

		\enr\setl \item[1.] $v(a \oid b)=0$, $v(a \noid b)=0$ (incomplete information), \item[2.]
		$v(a \oid b)=0$, $v(a \noid b)=1$ (only negative information), \item[3.] $v(a \oid b)=1$,
		$v(a \noid b)=0$ (only positive information), \item[4.] $v(a \oid b)=1$, $v(a \noid b)=1$
		(contradictory information), \eenr

	\mh and two scenarios of reliable information (i.e. when $v(\cons (a \oid b)) = 1$),

		\enr\setl

		\item[5.] $v(a \oid b)=1$, $v(a \noid b)=0$,

		\item[6.] $v(a \oid b)=0$, $v(a \noid b)=1$.

		\eenr

	\noi Scenarios (5) and (6) are classical, so $\langle \I(a), \I(b) \rangle$ must be either in the
	extension or in the anti-extension of $\oid$, not both.  The sensible point is how to express the
	non-classical scenarios (1) (no information at all) and (4) (conflicting information).  Let us
	see how they are represented in the semantics of \qletf.

	Let $a$, $b$, and $c$ be the names `Hesperus', `Phosphorus', and `Venus', and $P$ the predicate `is a
	planet'. Now consider the following hypothetical scenarios: \enr

		\item[(i)] $\{ \neg Pa, \neg Pb, Pc \}$.

		There is the information that neither Hesperus nor Phosphorus is a planet, and Venus is a
		planet.  Nothing is said about whether or not they are the same object. Thus, if $\I(a) =
		\bar{a}$ and $\I(b) = \bar{b}$, the pair $\langle \bar{a}, \bar{b} \rangle$ is not in
		$\oid_{-}$ and, of course, nor in $\oid_{+}$.

		\item[(ii)] $\{ a \noid b, \neg Pa, \neg Pb, Pc \}$.

		Scenario (ii) is like  to (i), except that we have the additional information that
		Hesperus and Phosphorus are not the same object. Thus, the pair $\langle \bar{a}, \bar{b}
		\rangle$ is in $\oid_{-}$.

		Regarding scenarios (i) and (ii), note that according to the semantic clauses of
		identity, $\langle \I(a),\I(b) \rangle $ belongs to the anti-extension $\oid_{-}$ if and
		only if $a \not\oid b$ holds. The anti-extension $\oid_{-}$ does not contain every pair
		$\langle \bar{a},\bar{b} \rangle$ such that $\I(a) \neq \I(b)$, but only the pairs
		$\langle \bar{a},\bar{b} \rangle$ such that $v(a\not\oid b)=1$.

		\item[(iii)] $\{ a \oid b, \neg Pa, \neg Pb, Pc \}$.

		Scenario (iii) is also like (i) except that now we have the additional  information that Hesperus and
		Phosphorus are indeed the same object. Thus, $\I(a) $ and $\I(b)$ have to have the same
		denotation, say, $\bar{a}$, and the pair $\langle \bar{a}, \bar{a} \rangle$ is, of course, in
		$\oid_{+}$ (as well for every other object in the domain).

		\item[(iv)] $\{ a \noid b, \neg Pa, \neg Pb , a \oid c, b \oid c, Pc \}$.

		This is a contradictory scenario that adds to scenario (ii) the information that both
		Hesperus and Phosphorus are in fact the planet Venus, and so the same object.

		Here we have contradictory information about their identity,  which means that $a$ and
		$b$ denote one and the same object, say $\bar{a}$, but also that these names denote
		different objects.  To express this scenario we make $\I(a) = \I(b)=\bar{a}$, and the
		pair $\langle \bar a, \bar a \rangle$ is in both $\oid_{+}$ and $\oid_{-}$.

		\eenr

	Two additional remarks are in order here. First, structures of \qletf\ may be thought of as
	representations of databases with positive and negative information that is marked as either reliable or
	unreliable. In this way, a structure is \textit{determined} by a configuration of a database. Scenario
	(i) above, for example, would correspond to the structure below: \enr \item[] $D = \{ \bar{a},
	\bar{b},\bar{c} \}$,

		$\I(a) = \bar{a} $, $ \I(b) = \bar{b}$, $ \I(c) = \bar{c}$

		$\I(P_+) = \{ \bar{c} \}$, $\I(P_-) = \{\bar{a}, \bar{b} \}$,

		$\I(\oid_+) = \{ \langle \bar{a}, \bar{a}\rangle \, \langle \bar{b}, \bar{b}\rangle, \langle \bar{c}, \bar{c}\rangle \} $,

		$\I(\oid_-) = \emptyset$.\footnote{Note that this structure makes $v(\incon (a \oid b))
		=1$, but leaves undetermined the values of $\con Pa$, $\con Pb$, and $\con Pc$, which
		have to be established by means of valuations.}

	\eenr

	\noi Indeed, the basic idea of the intuitive interpretation in terms of evidence/information is
	that sentences merely provide positive and negative information about objects and their
	properties, which of course does not mean that these sentences are true, nor that the putative
	objects they refer to exist. This brings us to the second remark. In the classical scenarios (5)
	and (6) we may assume that the names denote `real objects', that is, the sentences $a \oid b$ and
	$a \noid b$ talk about objects that indeed exist in the world. On the other hand, it might well
	be that in the non-classical scenarios the object that the interpretation assigns to a name does
	not exist. It is not difficult to imagine, for example, a database that contains information
	about an individual that, in fact, does not exist. The domain, in this case, reflects this
	situation and has an object that corresponds to that name, i.e.  that is the pseudo-denotation of
	the name.

	\end{rmrk}

\section{First-order \textit{FDE} and some of its extensions} \label{sec:qfdeetc}

Now that we have proven soundness and completeness theorems for \qletf, we will show in this section how
\qletf\ can be modified to yield first-order versions of $FDE$ and some of its extensions, namely,
Kleene's \textit{K3}, the logic of paradox \textit{LP}, and classical logic -- to be called here \qfde,
\qk, and \qlp, respectively. We start by presenting a natural deduction system and a corresponding
first-order valuation semantics for \qfde, which result from slight modifications of those for \qletf. We
shall also indicate how the completeness proof above can be adapted to the case of \qfde\ and, in
Section~\ref{sec:qfde-extensions}, to those of \qk\ and \qlp.

\subsection{First-order \textit{FDE}}

The logical vocabulary of \qfde\ is the same as that of \qletf\ except that $\con$ and $\incon$ are no
longer included in the set of logical primitives. We will continue to make use of first-order signatures
to specify the non-logical vocabulary of a first-order language, and adopt the same notational
conventions as before.

	\begin{df}\label{df:fdedns} The logic \qfde\ is obtained by dropping rules $EXP^{\con}$, $PEM^{\con}$,
	$Cons$, $Comp$, and $AV$ from \qletf\ (see Definition \ref{df:nds}).
	\end{df}

Notice that \qfde\ includes all $\con$- and $\incon$-free rules of \qletf, except for $AV$. As it turns
out, the absence of the $\con$ and $\incon$ allows to prove that any two alphabetically variant sentences
are deductively equivalent.  \black

Recall that while presenting the semantics of \qletf\ in Section \ref{sec:qletf} it was necessary to
supplement a structure \A\ with a valuation $v$ in order to ensure that all sentences in which $\con$ or
$\incon$ occur get assigned a semantic value -- since their values are not always determined by the
values of their subformulas. In \qfde, however, valuations are no longer necessary. Thus, although the
notion of a first-order structure remains the same as before (see Definition \ref{df:structure}),
interpretations, in the sense of Definition \ref{df:interpretation}, could be dispensed with in the case
of \qletf. Nonetheless, in order to preserve the notation used in the preceding sections, and to
demonstrate the generality of of the method of anti-extensions + valuations,
we shall continue to talk as if sentences get assigned one of the values $1$ or $0$ by a valuation in
\qfde, but this time to each structure \A\ there will correspond a \textit{single} valuation $v_{\A}$,
which is the valuation induced by \A.

	\begin{df}\label{df:fdevaluation} Let $\sig = \langle \mathcal{C}, \mathcal{P} \rangle$ be a
	signature and \A\ an \sig-structure. The mapping $v_{\A}: Sent(\sig_{\A}) \longrightarrow \{0,
	1\}$ is the \textit{valuation induced by} \A\ if it satisfes clauses (1)-(7) and (10)-(13) of
	Definition \ref{df:valuation} (where $v$ is replaced everywhere by $v_{\A}$).\\

	\noi Given $A \in Sent(\sig)$ and an \sig-structure \A, we shall say that $A$ \textit{holds} in
	\A\ ($\A \vDash A$) if and only if $v_{\A}(A) = 1$; and that $A$ is a \textit{semantic
	consequence of} $\Gamma$ in \qfde\ if and only if $\A \vDash A$ whenever $\A \vDash B$, for every
	$B \in \Gamma$.
	\end{df}

By suitable modifications of the definitions and the proofs of the results in Section
\ref{sec:soundnessandcompleteness}, it can be proven that \qfde\ is sound and complete with respect to
the class of all \qfde-structures. Specifically, the proof of (the \qfde-analogue of) Proposition
\ref{prop:substitution_02} (and so of Corollary \ref{cor:substitution_03}) is the same as before. Hence,
except for absense of proofs for the rules involving $\con$ and $\incon$, the soundness proof for \qfde\
remains the same as that for \qletf. As for completeness, there are no significant changes either. In
particular, the proof of (the \qfde-analogue of) Lemma \ref{lm:canonical} differs from the one presented
above only in that we are not required to show that $v_{\A}$, where \A\ is the canonical structure of a
regular Henkin set, satisfies clauses (8), (9), (14), and (15).

	\begin{rmrk}\label{rmrk:omori}\textit{On a constructive first-order \fde}

	\noi A sequent calculus for quantified \fde\ is found in Anderson and Belnap \cite{fde63}, and
	natural deduction systems in Priest \cite[pp.~331ff.]{priest2002} and Sano and Omori
	\cite[p.~463]{omori.qfde}.  The deductive systems in \cite{fde63} and \cite{priest2002} are
	equivalent to the one above, but that of \cite{omori.qfde}, as far as we can tell, is not. It
	seems to us not only that (i) there is a gap in the completeness proof of \cite{omori.qfde} but
	also that (ii) the natural deduction system presented therein is in fact incomplete.
	Specifically, it cannot prove:
		\begin{align}
		\tag{1}\label{tag:sano&omori_02} \forall x(B \lor A) \vdash B \lor \forall x A,
		\end{align}

	\noi where $x$ is not free in $B$.
	Items (i) and (ii) bring to light some interesting points about a constructive first-order
	extension of \fde.

Let us take a look at (i). In \cite{omori.qfde} Sano and Omori introduce some formal systems that extend
the first-order version of $FDE$ (\qfde), which they call $BD$ logic. They present a natural deduction
system for \qfde, where the introduction rule for the universal quantifier is

		{\small \begin{center}

			\bottomAlignProof \AxiomC{$A(c/x)$}
			\RightLabel{$\forall I'$} \UnaryInfC{$\forall xA$} \DisplayProof

		\end{center}}

	\noi instead of 

\begin{center} \small
	\bottomAlignProof
					\AxiomC{$B \lor A(c/x)$}
				\RightLabel{$\forall I$}
				\UnaryInfC{$B \lor \forall xA$}
			\DisplayProof
	\end{center}

\noi (with the usual restrictions). The semantics is similar to the one above
	except for the use of relations instead of functions. \citet[Sect. 4]{omori.qfde} offer a general method
	for proving the completeness of $BD$ and some of its extensions, but the proof of the result corresponding
	to Lemma \ref{lm:lindenbaum} above  \cite[Lemma~4.2]{omori.qfde} is not entirely convincing. They adopt the
	following notation \citep[p.~464]{omori.qfde}:
		\begin{align} \tag{2}\label{tag:sano&omori_01} \Gamma \vdash_\mathcal{R} \Pi \mbox{ iff for some
		finite subset } \{ A_{1},...,A_{j} \} \mbox{ of } \Pi, \Gamma \vdash_\mathcal{R} A_{1} \lor ...
		\lor A_{j} \end{align}

	\noi for any calculus $\mathcal{R}$ such that $BD \subseteq \mathcal{R}$. They then proceed with a
	Lindenbaum construction, defining a sequence $\langle \Gamma_{n}, \Pi_{n} \rangle$, and claim, without
	presenting a proof, that for every $n$, $\Gamma_{n} \nvdash_{\mathcal{R}} \Pi_{n}$.\footnote{In \cite[p. 466]{omori.qfde} we just read that ``[for every $n$] it is easy to see that
	$\Gamma_n \nvdash_\mathcal{R} \Pi_n $''.} The sensible point of their strategy is to show that
	$\Gamma_{n+1} \nvdash_{\mathcal{R}} \Pi_{n+1}$ when $\Gamma_{n}, \forall xB \vdash \Pi_{n}$. In this case,
	$\Gamma_{n+1} = \Gamma_{n}$ and $\Pi_{n+1} = \Pi_{n} \cup \{\forall xB, B(c/x)\}$, where $c$ is a new
	constant. If we were to fill in the gaps in their proof, we could assume that (a) $\Gamma_{n}
	\vdash_{\mathcal{R}} \Pi_{n} \cup \{\forall xB, B(c/x)\}$ to obtain (b) $\Gamma_{n} \vdash_{\mathcal{R}}
	\Pi_{n}$, which yields a contradiction with the induction hypothesis. However, there is no obvious way to
	get from (a) to (b) by applying rule $\forall I'$ instead of $\forall I$.

  Concerning (ii), let us call \qfde$'$ the first-order extension of \fde\ with rule $\forall I'$ instead of $\forall I$.
	As far as we can see,  \qfde$'$ is incomplete with respect to the standard
	semantics for quantified \fde, an issue which is closely related to the validity of
	(\ref{tag:sano&omori_02}). A sketch of a proof that \qfde$'$ dos not prove (\ref{tag:sano&omori_02}) is as follows. 
	{Let us call \qfde$_G'$ the $\to$-free   fragment of L\'opez-Escobar's refutability calculus  	\cite{escobar}. It is straightforward to prove that if $\Gamma\vdash A$ does not hold  in \qfde$_G'$, then it does not hold in 
	\qfde$'$.  Now, define a notion of generalized subformula in such a way that $\neg A$ and $\neg B$ are generalized subformulas 
	of $\neg(A \land B)$, $\neg A(t/x)$ of $\neg\forall x A $, and so on. Since cut-elimination holds for \qfde$_G'$, it is easy to see that all formulas in a cut-free derivation in \qfde$_G'$  are generalized subformulas 
	of the endsequent of the derivation. So, if (\ref{tag:sano&omori_02}) were valid in \qfde$_G'$, 
	it would be provable with the positive rules only, but every proof-search ends with a topsequent which is not an axiom.}  It is also 
	worth 	noting that if (\ref{tag:sano&omori_02}) were provable with the positive fragment 
	of  \qfde$_G'$, it would be provable in intuitionistic logic, but it is not. 
		Indeed, the natural way of constructively extending \fde\ to a first-order logic is given by \qfde$'$.
Notice, besides, that contrary to \qfde$'$, the rules of \qfde\ are not 	harmonious precisely because the elimination rule $\forall E$ cannot be `read off' from $\forall I$.  
A corresponding adequate semantics for \qfde$'$ would not be given by its standard
	semantics, found here and in \cite{priest2002,omori.qfde}.  

	\end{rmrk}

\subsection{On Some Extensions of \qfde}\label{sec:qfde-extensions}

It should be no surprise by now that the definitions and results presented in sections \ref{sec:qletf} and
\ref{sec:soundnessandcompleteness} can also be straightforwardly modified to yield sound and complete natural
deduction systems for some well-known extensions of \qfde, namely, the first-order versions of the logic of
paradox \textit{LP} \cite{priest.lp} and \cite[Ch.~5]{priest.ic}, Kleene's \textit{K3} \cite{kleene}, and even
classical logic. One has to simply add either excluded middle or explosion (or both) to the rules of \qletf. As
for the semantics, it suffices to impose some further conditions on the relations between the extensions and
anti-extensions of predicate letters.

	\begin{df}\label{df:qfde-extensions-nds} Consider the following two rules:

		{\small \begin{center}

			\bottomAlignProof
					\AxiomC{}
				\RightLabel{$PEM$}
				\UnaryInfC{$A \lor \neg A$}
			\DisplayProof
		\qquad
			\bottomAlignProof
					\AxiomC{$A$}
					\AxiomC{$\neg A$}
				\RightLabel{$EXP$}
				\BinaryInfC{$B$}
			\DisplayProof

		\end{center}}
		\begin{itemize}\setl

		\item[1.] \qlp\ results from adding $PEM$ to the rules of \qfde;

		\item[2.] \qk\ results from adding $EXP$ to the rules of \qfde;

		\item[3.] \qcl\ results from adding both $PEM$ and $EXP$ to the rules of \qfde.

		\end{itemize}

	\end{df}

	\begin{df}\label{df:qfde-extensions-interpretation} Let \sig\ be a signature and let $\A = \langle \D, \I
	\rangle$ be a \qfde-structure. Then:

		\begin{itemize}\setl

		\item[E1.] \A\ is a \qlp-\textit{structure} if and only if $P_{+}^{\A} \cup P_{-}^{\A} = \D^{n}$,
		for every $n$-ary predicate letter $P$ of \sig;

		\item[E2.] \A\ is a \qk-\textit{structure} if and only if $P_{+}^{\A} \cap P_{-}^{\A} =
		\emptyset$, for every predicate letter $P$ of \sig;

		\item[E3.] \A\ is a \qcl-\textit{strucure} if and only if it satifies both (E1) and (E2).
		\end{itemize}

	\noi If $\mathcal{L}$ is one of \qlp, \qk, or \qcl, then an $\mathcal{L}$-\textit{interpretation} is a
	pair $\langle \A, v_{\A} \rangle$, where \A\ is an $\mathcal{L}$-struture and $v_{\A}$ is the
	\textit{valuation induced by} \A\ -- i.e., mapping from $Sent(\sig)$ to $\{1,0\}$ that satisfies clauses
	(1)-(7) and (10)-(13) of Definition \ref{df:valuation}.

	\end{df}

This way of presenting classical predicate logic is 
unusual and filled with redundancies --
\qcl\ could be more simply described as the logic resulting from adding $EXP$ and $PEM$ to the positive fragment
of \qfde\ with the usual introduction rule for $\forall$ (i.e., $\forall I'$). Moreover, the semantics, and so the
metatheoretical results used above for proving completeness, could also be simplified. As with \qfde, valuations
are not strictly necessary, since structures alone suffice to determine the semantic values of all formulas -- and
in the case of \qcl\ anti-extensions are not   required. Our point in adapting the semantics of \qletf\ to
\qfde, \qk, \qlp, and \qcl, however,  is to show that
the  method of anti-extensions + valuations can be easily applied to a number of non-classical logics, and even to classical logic.

\begin{lm}\label{lm:qfde-extensions-props} Let \A\ be a \qfde-structure. Then:

		\begin{itemize}\setl

		\item[1.] If \A\ is a \qlp\ or a \qcl-structure, then $v_{\A}(A) = 1$ or $v_{\A}(\neg A) = 1$, for
		every $A \in Sent(\sig_{\A})$;

		\item[2.] If \A\ is a \qk- or a \qcl-structure, then $v_{\A}(A) \neq 1$ or $v_{\A}(\neg A) \neq
		1$, for every $A \in Sent(\sig_{\A})$.

		\end{itemize}

	\end{lm}

	\begin{pr} Both results follow by straightforward inductions on the complexity of $A$. The atomic cases of
	the proofs of (1) and (2) depend respectively on conditions (E1) and (E3), and on (E2) and (E3) of
	Definition \ref{df:qfde-extensions-interpretation}.

	\end{pr}

	\begin{tr}\label{tr:qfde-extensions-completeness} Let $\mathcal{L}$ be one of \qlp, \qk, or \qcl. Then,
	$\mathcal{L}$ is sound and complete with respect to the class of all $\mathcal{L}$-structures.

	\end{tr}

	\begin{pr} Lemma \ref{lm:qfde-extensions-props}(1) is all we need to prove that rule $PEM$ is valid in
	both \qlp\ and \qcl, while the validity of $EXP$ in \qk\ and \qcl\ follows from
	\ref{lm:qfde-extensions-props}(2). By adapting the proof of the soundness of \qletf\ (Theorem
	\ref{tr:soundness}), it can then be easily proven that all three systems are sound with respect to the
	class of corresponding structures. Proving their completeness is equally straightforward: it suffices to
	ensure that the canonical structure defined in the proof of Lemma \ref{lm:canonical} satisfies the
	corresponding restriction in Definition \ref{df:qfde-extensions-interpretation}. But this is an immediate
	consequence of the presence, in each case, of either $PEM$ or $EXP$ (or both) in the deductive
	system\footnote{In the case of \qlp, for example, given a regular Henkin set $\Delta$, it follows by $PEM$
	that $\Delta \vdash_{\qlp} P(c_{1},\dots,c_{n}) \lor \neg P(c_{1},\dots,c_{n})$, and so that either
	$\Delta \vdash_{\qlp} P(c_{1},\dots,c_{n})$ or $\Delta \vdash_{\qlp} \neg P(c_{1},\dots,c_{n})$.
	Therefore, if $\A = \langle \D, \I \rangle$ is $\Delta$'s canonical structure, then for every
	$a_{1},\dots,a_{n} \in \D$, either $\langle a_{1},\dots,a_{n} \rangle \in P_{+}^{\A}$ or $\langle
	a_{1},\dots, a_{n} \rangle \in P_{-}^{\A}$.}.
	\end{pr}

\section{Final Remarks} \label{sec:final}

The so called Suszko's thesis \cite{suszko77} asserts that every Tarskian and structural logic admits of
a two-valued semantics. A proof of this result for sentential logics can be found in
\cite[pp.  72-73]{malinowski}, and is in fact very simple. Given a (possibly infinite) multi-valued semantics for
a Tarskian and structural logic $\mathcal{L}$, a two-valued semantics for $\mathcal{L}$ is defined as follows: if
a formula $A$ receives a designated valued in a multi-valued interpretation $\I$, the value 1 is assigned to $A$
in a two-valued interpretation $\I'$, otherwise $A$ is assigned the value 0 in $\I'$. Semantic consequence is then
defined as preservation of the value 1, instead of preservation of a designated value.

An analogous result has been obtained  by Loparic and da Costa in \cite[pp.~121-122]{costa_loparic},
where they present a general notion of valuation semantics. Given a consequence relation $\vdash$ and a language
$\mathsf{L}$, a function $e: \mathsf{L} \longrightarrow \{0,1\}$ is an \textit{evaluation} if $e$ satisfies the
following clauses:
	\begin{description}\setl
	\item (i) If $A$ is an axiom, then $e(A) = 1$;
	\item (ii) If $e$ assigns the value $1$ to all the premises of an application of an inference rule, then 	it also assigns $1$ to its conclusion;
	\item (iii) For some formula $A$, $e(A) = 0$.
	\end{description}
\noi It is also necessary  that a Lindenbaum construction can be carried out for $\vdash$, which requires that
 $\vdash $ has to be Tarskian and compact.  Let a set $\Delta$ be $A$-\textit{saturated} when $\Delta \nvdash A$ and for every $B \notin \Delta$, $\Delta \cup \{ B\} \vdash A$. Now, assuming that $\Gamma \nvdash A$:
	\begin{description}\setl
	\item (iv) There is an $A$-saturated set $\Delta$, such that $\Gamma \subseteq \Delta$;
	\item (v) $\Delta \vdash B$ iff $B \in \Delta$;
	\item (vi) The characteristic function $c$ of $\Delta$ is an evaluation.
	\end{description}
\noi Since (iv) and (v) are immediate consequence of the Lindenbaum construction, it suffices to prove (vi).
Clearly, $c$ satisfies (i) and (iii) above. As for (ii), suppose $c$ assigns the value 1 to the premises of a
derivation $\Delta_{0} \vdash B$, $\Delta_{0} \subseteq \Delta$. Since, by (v), $B \in \Delta$, it then follows
that $c(B) = 1$.

Now, define a \textit{valuation} as an evaluation that is the characteristic function of some $A$-saturated set.
The collection of all valuations so defined is an adequate \textit{valuation semantics} for $\vdash$.  Soundness
follows from the definition of evaluations (the set of valuations is a proper subset of the set of evaluations),
and completeness from the fact that $c$ assigns $1$ to all the sentences of $\Gamma$, while assigning $0$ to
$A$.\footnote{Note that the set of all evaluations for a given consequence
relation $\vdash$ does not suffice for providing a semantics. Consider e.g.  the semantics of classical logic,
which is a special case of a valuation semantics, and let $T$ be the set of all classical theorems. The
characteristic function $c$ of $T$ is an evaluation, but for all atoms $p$, neither $p$ nor $\neg p$ is in $T$,
so $c(p)=0$ and $c(\neg p)=0$, even
though $c(p \lor \neg p)=1$.
It is also worth noting that the notion of an $A$-saturated set provides a
method for proving completeness for any logic for which a Lindenbaum construction can be carried out. We thank
 Andrea Loparic for some conversations that clarified the general notion of valuation semantics.}
 Apparently, provided appropriate conditions for the construction of an $A$-saturated set, this result could
be extended to first-order logics as well.

Valuation semantics were proposed by Loparic, Alves and da Costa for the paraconsistent
logics of da Costa's $Cn$
hierarchy \cite{costa.alves,loparic1986,loparic.alves}, which are `ancestors' of the logics of formal
inconsistency and logics of evidence and truth.
The problem they had at hand was to provide semantics for  paraconsistent
logics that are not finitely-valued.
They then came up with the idea of generalizing classical
two-valued semantics in such a way that
the axioms and rules were `mirrored' by the semantic clauses
in terms of $0$s and $1$s.
The value $0$ assigned to a
formula $A$ can be read as `$A$ does not hold' and $1$ as `$A$ holds' --
note that this is the basic idea
of the general notion of valuation as defined above.
Later, valuation semantics were proposed for several
non-classical sentential logics, including minimal and intuitionistic logic,
\fde, Nelson's \nel, and logics of
formal inconsistency and undeterminedness
\cite{letj, carcoma2007,loparic2010,costa_loparic,letf}.
First-order valuation semantics were also proposed for da Costa's
quantified $Cn$ hierarchy  \cite{bueno.costa.krause},
and for some logics of formal inconsistency \cite{carco.book,qmbc}.

As we have seen above (Definition \ref{df:interpretation} and
Remark \ref{rmrk:qletf.inte}), when a non-deterministic
semantics is extended to first-order  the crucial point is
how to handle its extended non-deterministic character.
The results presented here
suggest {anti-extensions + valuations} as
a general method for providing first-order valuation semantics for
non-classical logics. Of course, these tools can sometimes be simplified, as we have just seen in Section \ref{sec:qfdeetc}.  Valuations can be dispensed with in the case of \qfde, \qlp, and \qk, but are indispensable in
the case of \qletf, along with several other logics of formal inconsistency (e.g. \qmbc\ \cite{qmbc}). Classical
logic is a limiting case, since standard Tarskian structures (where anti-extensions are just the complement of extensions)
are enough to provide an adequate semantics.
In all these cases, however, the semantics
are nothing but special cases of the general method described here.

\bibliographystyle{plainnat}


\end{document}